\documentclass[preprint,3p]{elsarticle}

\usepackage{amssymb}
\usepackage{amsmath}
\usepackage{bm}

\usepackage{xcolor}
\usepackage{hyperref}
\usepackage{xspace}
\usepackage{subcaption}
\usepackage{tikz}

\biboptions{sort&compress}

\newcommand{\argmin}{\mathop{\mathrm{argmin}}}

\newcommand{\ddfrac}[2]{\frac{\partial #2}{\partial #1}}

\newcommand{\rhou}{  \rho u }

\newcommand{\rhoe}{  \rho e }

\newcommand{\ddx}{\frac{\partial}{\partial x}}

\newcommand{\ddt}{\frac{\partial}{\partial t}}
\newcommand{\ddr}{\frac{\partial}{\partial r}}

\newcommand{\cEDNN}{C-EDNN\xspace}
\newcommand{\dEDNN}{D-EDNN\xspace}
\newcommand{\mEDNN}{Multi-EDNN\xspace}

\newcommand{\NUMDIM}{S}
\newcommand{\NUMSTATE}{K}
\newcommand{\NUMWEIGHTS}{N}
\newcommand{\NUMX}{M}
\newcommand{\NUMELES}{E}
\newcommand{\NUMBOUND}{B}

\newcommand{\R}{\mathbb{R}}
\newcommand{\Rplus}{\mathbb{R}^+}
\newcommand{\Rdim}{\mathbb{R}^\NUMDIM}
\newcommand{\Rstate}{\mathbb{R}^{\NUMSTATE}}

\newcommand{\xdom}{\bm{x}}
\newcommand{\xbou}{\bm{\chi}}
\newcommand{\weights}{\bm{W}}
\newcommand{\soln}{\bm{q}}
\newcommand{\aux}{\bm{a}}
\newcommand{\flux}{\bm{f}}

\newcommand{\dis}{D}
\newcommand{\cor}{C}
\newcommand{\operator}{\bm{\mathcal{N}}}
\newcommand{\domain}{\bm{\Omega}}
\newcommand{\boundary}{\bm{\partial\Omega}}
\newcommand{\fluxfunc}{\bm{\mathcal{F}}}
\newcommand{\solnSolver}{\bm{\mathcal{Q}}^\star}
\newcommand{\fluxSolver}{\fluxfunc^{\perp\star}}

\newcommand{\jacobian}{\bm{J}}
\newcommand{\rhs}{\bm{N}}
\newcommand{\weightUpdate}{\bm{Z}}

\newcommand{\solnCorrectFunc}{g}
\newcommand{\fluxCorrectFunc}{\bm{h}}

\newcommand{\bcWeightSolid}{c_{s}}

\newcommand{\adaptHeight}{h}

\newcommand{\adaptArea}{A}
\newcommand{\adaptParam}{\beta}

\journal{Journal of Computational Physics}

\begin{document}

\begin{frontmatter}



\title{Multi evolutional deep neural networks (Multi-EDNN)}

% \author[label1,label2]{Hadden Kim, Tamer Zaki}
% \affiliation[label1]{organization={},

% \affiliation[label2]{organization={},

\author[ams]{Hadden Kim}
\author[ams,me]{Tamer A.\,Zaki\corref{cor1}}

\affiliation[ams]{organization={Department of Applied Mathematics \& Statistics, Johns Hopkins University},%Department and Organization
            addressline={3400 N.\,Charles St.}, 
            city={Baltimore},
            postcode={21218}, 
            state={MD},
            country={USA}}
\affiliation[me]{organization={Department of Mechanical Engineering, Johns Hopkins University},%Department and Organization
            addressline={3400 N.\,Charles St.}, 
            city={Baltimore},
            postcode={21218}, 
            state={MD},
            country={USA}}
\cortext[cor1]{t.zaki@jhu.edu}

\begin{abstract}
Evolutional deep neural networks (EDNN) solve partial differential equations (PDEs) by marching the network representation of the solution fields, using the governing equations.  Use of a single network to solve coupled PDEs on large domains requires a large number of network parameters and incurs a significant computational cost.  We introduce coupled EDNN (\cEDNN) to solve systems of PDEs by using independent networks for each state variable, which are only coupled through the governing equations.  We also introduce distributed EDNN (\dEDNN) by spatially partitioning the global domain into several elements and assigning individual EDNNs to each element to solve the local evolution of the PDE.  The networks then exchange the solution and fluxes at their interfaces, similar to flux-reconstruction methods, and ensure that the PDE dynamics are accurately preserved between neighboring elements.  Together \cEDNN and \dEDNN form the general class of \mEDNN methods.   We demonstrate these methods with aid of canonical problems including linear advection, the heat equation, and the compressible Navier-Stokes equations in Couette and Taylor-Green flows.  

\end{abstract}


% \begin{highlights}
% \item Coupled evolutional deep neural networks (C-EDNN) for solving systems of PDEs.
% \item Distributed EDNN (D-EDNN) with interface treatments for PDEs on large domains.
% \item Coupled, distributed networks (Multi-EDNN) demonstrated for compressible Navier-Stokes.
% \end{highlights}

\end{frontmatter}


%%%%%%%%%%%%%%%%%%%%%%%%%%%%%%%%%%%%%%%%%%%%%%%%%%%%%%%%%%%%%
%   Main Text
%%%%%%%%%%%%%%%%%%%%%%%%%%%%%%%%%%%%%%%%%%%%%%%%%%%%%%%%%%%%%
\section{Introduction}
\label{sec:intro}

Use of machine learning (ML) in the field of partial differential equations (PDEs) and physics-based predictions has seen several exciting developments.  
Some efforts have focused on learning from data and introducing physics constraints \cite{Raissi_19, Jagtap_20, Patricio_23}, and others have been developed to learn operators \cite{Lu_21, Li_20, Cai_21}.
The particular area of interest for the present work is use of machine learning to numerically solve partial differential equations (PDEs) and, in particular, to forecast the evolution of dynamical systems. 
This idea was pursued by \citet{Du_ednn}, who introduced evolutional deep neural networks (EDNN, pronounced ``Eden'').  They used the expressivity of neural networks to represent the solution of the PDE in all but the direction of propagation, and then adopted the age-old marching approach to evolve the network parameters according to the governing PDEs.  In doing so, the evoltion of the network becomes predictive, providing a forecast for any horizon of interest.  
When the evolution is in time, the EDNN weights can be viewed as traversing the space of network parameters through paths directed by the governing PDE (see figure \ref{fig:ednn}).  In the present work, rather than deploying a single network to solve the system of PDEs on the entire domain, we introduce a decomposition of the problem, both in state space and in the spatial domain, and utilize multiple concurrent EDNNs to solve the PDEs.  This contribution is essential to enable future use of ML-based techniques such as EDNN for the solution of complex and large-scale problems. 

Since the introduction of EDNN, a growing community has extended its capabilities.  Here we summarize a selection of related work that are most relevant to the present effort.  
\citet{Anderson_24} provided the insightful interpretation of EDNN as a reduced-order nonlinear solution.  The output layer of the deep neural network can be viewed as the combination of time-dependent modes. Unlike in classical, linear reduced-order models, these modes shape-morph with time and are allowed to depend nonlinearly on their parameters. 
\citet{Bruna_23} interpretted EDNN as an active learning algorithm that can generate its own training data.   When starting from random parameters, neural networks may need a large number of training epochs to converge to the solution.  If, however, the network parameters are initialized close to the target, a single training step with very little training data may be sufficient.  In many PDE problems, the solution at one time is similar to a previous time. EDNN capitalizes on this, and its parameter evolution is viewed as training the parameters sequentially one time-step to the next with training data generated by EDNN itself.

As a meshless method, EDNN can be sampled anywhere in the spatial domain, removing the need for complex remeshing algorithms.  Researchers adopted adaptive sampling strategies to both reduce the number of points needed and to improve  accuracy in the network evolution.  \citet{Bruna_23} sampled the spatial domain with time-dependent, adaptive measures estimated by Guassian mixture models.
\citet{Wen_24} continued in this direction and formulated the adaptive sampling as the dynamics of particles within the spatial domain.  At each time-step, their method uses the network prediction to determine the particle trajectories, then samples at the displaced particle positions. 
\citet{Kast_24} first collected the network prediction over a high-resolution finite-element mesh and formulated a probability distribution based on the magnitude of the PDE dynamics at the mesh points.  Then they used only a subsampling of these points in the evolution of EDNN.

Several advancements have been made to enable EDNN to solve PDEs with boundary conditions. 
\citet{Du_ednn} utilized an input feature layer to enforce periodicity. 
They also enforced Dirichlet boundary conditions by negating the boundary predictions of the network.
\citet{Kast_24} used a positional embedding (feature) layer to enforce both Dirichlet and Neumann boundary conditions.
With pre-computed embeddings, the EDNN prediction is restricted to the solution subspace that respects the boundary conditions of the problem.
\citet{Chen_23} utilized a simple and flexible approach:
by adding penalty terms to the EDNN algorithm, the network evolution can be driven to agree with the boundary conditions.

Other efforts have been directed to solution method itself.  
\citet{Anderson_24} adopted a constrained optimization to ensure that the EDNN evolution respects conservation laws.
Additionally, they introduced Tikhonov penalization to regularize the optimization objective function, which results in computational speedup. 
\citet{Chen_23} investigated a variety of time integrators, such as an energy-preserving midpoint method, implicit Euler, and Chorin-style operator splitting.
\citet{Finzi_23} transformed the spatial input into a scaled sinusoidal and applied Nystr\"om preconditioning.
\citet{Kao_24} utilized tensor neural network architectures and explored updating only a random subset of the network parameters each timestep.

All these related works share the common strategy of using a single neural network to represent the solution of interest over the entire spatial domain.
Relying on only one EDNN, however, has its drawbacks.
First, it is challenging to use a single network to predict the evolution of a vector-valued solution. 
The potentially different functional forms of the different states and the differences in their dynamical characteristics require a large number of parameters to accurately predict. 
Second, a large number of spatial samples is required to accurately evolve the parameters of the network for these problems. 
Both drawbacks contribute to the computational cost.
While parallel computations of a single, large deep neural network \cite{Ben_19} continue to garner attention, we adopt a different approach that attempts to reduce the computational cost and that is naturally amenable to parallel execution. 

Our strategy is inspired by classical finite element methods (FEM).  
FEMs partition the spatial domain into elements and approximate the solution in each element with the linear combination of a finite set of basis functions, such as polynomials or trigonometric polynomials.
For transient problems, time-dependent coefficients of the selected basis form a system of ordinary differential equations (ODEs) which can be marched with a choice of time integrator. 
Typically, a very small set of basis functions is selected, limiting expressivity.
To recover accurate solutions, the domain is partitioned into a very large number of elements, e.g.\,millions or billions in the context of simulations of complex  flows \cite{Merzari_23}.
Different techniques to ensure that the solution behaves correctly between elements have spawned numerous numerical methods.
For example, conforming finite element methods \cite{becker_fem} require the solution approximation be continuous everywhere, particularly on the element interfaces. 
Discontinuous Galerkin methods \cite{Cockburn_DG} allow the function to have jump discontinuities and utilize a Riemann solver to compute common numerical fluxes at element interfaces. 
Galerkin methods, however, seek weak solutions to the PDE which rely on a choice of the test function space. 
It is unclear which set of test functions would be appropriate for a deep neural network solution approximation. 
The flux-reconstruction family of methods \cite{Huynh_fr, Vincent_fr, Castonguay_13} lifts this requirement and instead formulates a correction to the local flux approximation of each element using data from neighboring elements. 
By reconstructing a flux whose normal components are continuous between elements, flux-reconstruction methods ensure global accuracy. 
The evolution of each element is formulated as a local procedure needing only minimal data synchronization at the element interfaces. 
As such, the element-specific parameters are nearly decoupled, which enables parallel and distributed computing \cite{Witherden_pyfr}.

In the present work, we develop a framework that is inspired by classical FEM to solve PDEs using coupled and distributed neural networks, rather than a single large network. We thus remove the limitation of using just one network to approximate the possibly vector-valued solution over the entire spatial domain.
First, we utilize coupled deep evolutional neural networks (\cEDNN) to represent each component of the state in a system of PDEs. 
Next, we partition the spatial domain and employ, for each subdomain-specific, distributed evolutional deep neural networks (\dEDNN) to collectively solve for the evolution of the PDE.
With these two strategies, we develop a framework of coordinating multiple EDNNs and call this method \mEDNN.
Section \ref{sec:prelim} gives a brief review of the original EDNN method.
In section \ref{sec:mednn}, we describe in detail the \mEDNN methodology, starting with coupling EDNNs for a system of PDEs and then we explain how we connect multiple EDNNs with domain decomposition. 
%  We also provide the construction for a new class of correction functions and our boundary condition procedure.
Section \ref{sec:experiments} presents the results of using \mEDNN to predict the evolution of the solution of a variety of PDE problems.
We summarize our conclusions in section \ref{sec:conclusion}.

\section{Problem setup and preliminaries}
\label{sec:prelim}

In this section, we briefly review the evolutional deep neural network (EDNN), first introduced by \citet{Du_ednn}.  
Consider a time-dependent, non-linear, partial differential equation
\begin{equation}
    \label{eq:pde}
    \ddfrac{t}{\soln} (\xdom, t) = \operator (\xdom, t, \soln)
\end{equation}
where the solution $\soln:\domain\times\Rplus \rightarrow \Rstate$ is a vector function on both the spatial domain $\domain \subset \Rdim$ and time $t \geq 0$, 
and $\operator$ is a differential operator. 
A time-dependent neural network is used to approximate the solution.
The inputs of the network are spatial coordinates, and the outputs are the solution approximated at the input coordinates.
By fixing the network topology, such as activation functions, number of layers, and number of neurons, the network is parameterized by its kernel weights and biases which are collectively denoted $\weights$.
We consider the network parameters as functions in time $\bm{W}:\Rplus \rightarrow \R^{\NUMWEIGHTS}$.
As such, the network can now be seen as the ansatz, 
\begin{equation}
    \label{eq:ansatz}
    \hat{\soln} \left(\xdom, \weights(t) \right).
\end{equation}
Using the chain rule, we connect the governing PDE to the time-evolution of the network parameters,  
\begin{equation}
    \label{eq:chain}
    \ddfrac{\weights}{\hat{\soln}} \frac{d\weights}{dt}
    % = \ddfrac{t}{\hat{\soln}} 
    = \operator (\xdom, t, \hat{\soln}).
\end{equation}
We can therefore express the time-evolution of the network parameters according to, 
\begin{equation}
    \label{eq:optimization}
    \hat{\bm{\zeta}}(t) = \argmin_{\bm{\zeta}} 
    \int \sum_{k=1}^{\NUMSTATE} \left\vert \ddfrac{\weights}{\hat{\soln}_k} \bm{\zeta} - \operator_k  \right\vert^2 d \domain.  
\end{equation}
Using a finite number $\NUMX$ spatial samples, 
we convert \eqref{eq:chain} into the linear system
\begin{equation}
\label{eq:linear}
    \jacobian \weightUpdate = \rhs,
\end{equation}
and the network evolution \eqref{eq:optimization} into
\begin{equation}    
    \label{eq:optimization_linear}
    \hat{\weightUpdate} = \argmin_{\weightUpdate} \lVert \jacobian \weightUpdate - \rhs \rVert_F^2
\end{equation}
where $\lVert \cdot \rVert_F$ is the Frobenius norm, $J_{k,i,j}$ is the network gradient, $Z_j$ is the weight update, and $N_{k,i}$ is the PDE operator.  The subscript $k$ identifies a solution variable within the state vector $\soln$, subscript $i$ is the spatial sampling point, and subscript $j$ is the network parameter.
Both $\jacobian$ and $\rhs$ are readily available by evaluating the neural network at a set of sampling points and using automatic differentiation to obtain any necessary derivatives.  By discretizing time and using classic time integrators, such as fourth-order Runge-Kutta (RK4), EDNN can now be seen as a time marching algorithm.  We sequentially evolve the weights from one time instant to the next using the trajectory $\hat{\weightUpdate}$.

\begin{figure}
    \centering
    \begin{subfigure}[b]{0.48\textwidth}
        \centering
        \includegraphics[width=\textwidth]{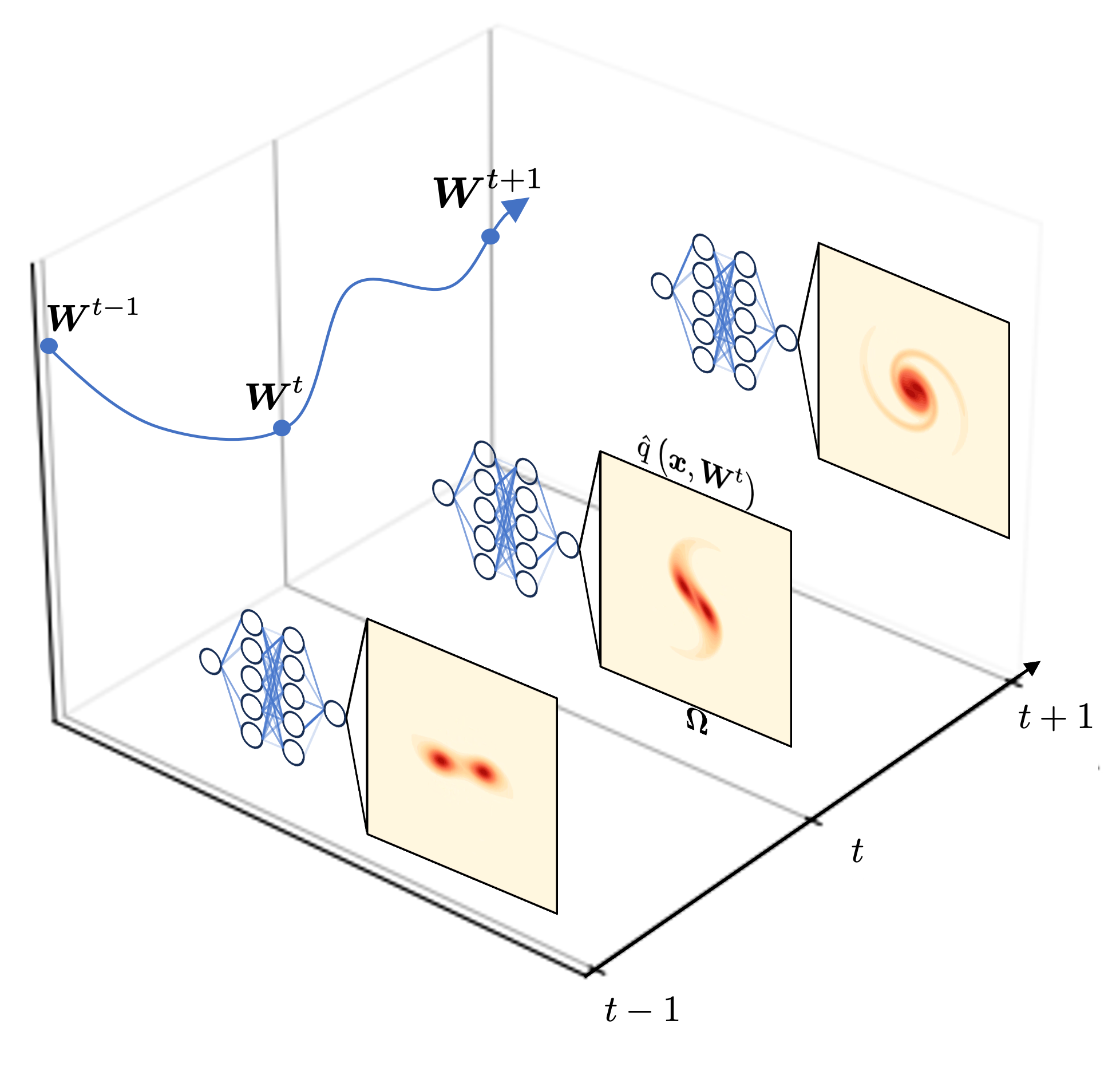}
        \caption{}
        \label{fig:ednn}
    \end{subfigure}
    \hspace{0.05\textwidth}
    \begin{subfigure}[b]{0.27\textwidth}
        \centering
        \includegraphics[width=\textwidth]{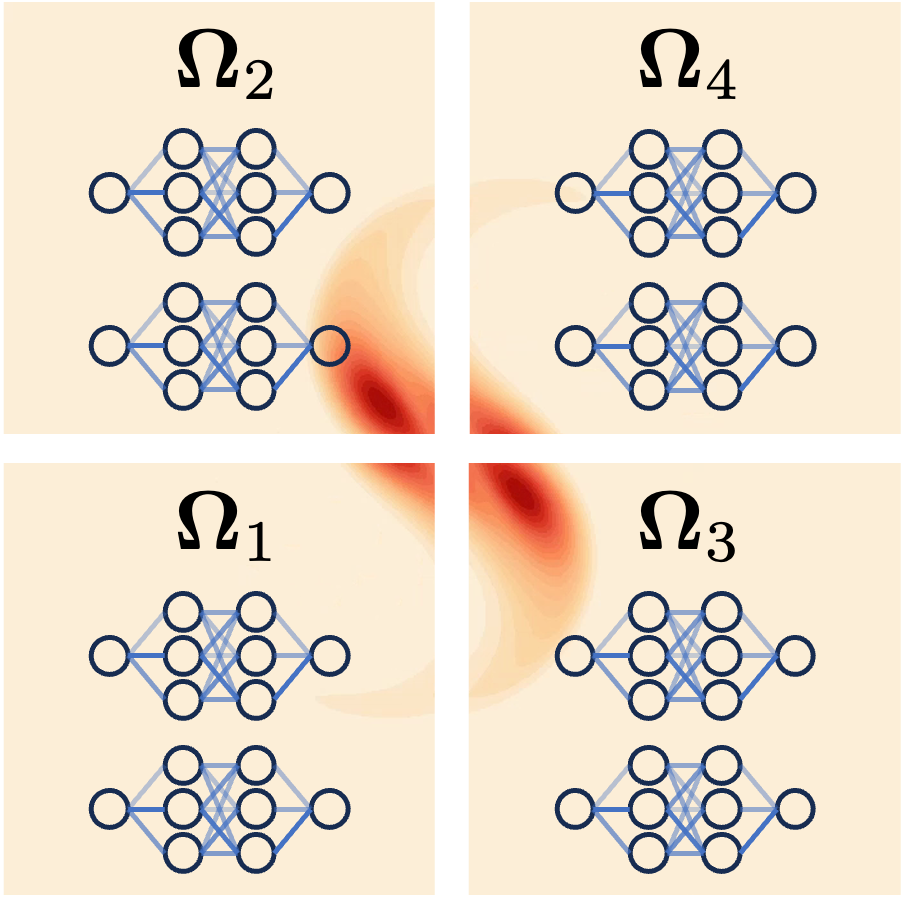}
        \vspace*{28pt}
        \caption{}
        \label{fig:mednn}
    \end{subfigure}
    \caption{Schematic of (a) original EDNN solving a system of PDEs on the entire spatial domain
    and (b) \mEDNN solving coupled differential equations using distributed networks on four subdomains.  In \mEDNN, each subdomain has two coupled networks that solve the system of equations locally, and that exchange information at the sub-domain interfaces with their neighbors. }
\end{figure}

A schematic of EDNN evolution over three time instances is shown in figure \ref{fig:ednn}.  At the first instance, the entire solution field is represented by the network $\hat{\soln} \left(\xdom, \weights^{t-1} \right)$. The weights are then evolved using the above governing equations to predict the next instance, $\hat{\soln} \left(\xdom, \weights^t \right)$.  The process is then repeated to predict $\hat{\soln} \left(\xdom, \weights^{t+1} \right)$ and can be continued for any time horizon of interest.

\section{Formulation of Multi Evolutional Deep Neural Networks (\mEDNN)}
\label{sec:mednn}

As originally introduced, EDNN utilizes a single deep neural network to represent the full state $\soln$ over the entire spatial domain. 
Problems that involve coupled fields, large domains, local variations in the behaviour of the solution, etc., require a large number of network parameters to accurately represent the solution and a large number of sampling points.
These two requirements contribute to the computational cost and memory requirements. 
In particular, the computational bottleneck of EDNN lies in the solving  equation \eqref{eq:optimization_linear}.

In order to address these challenges, we split the PDE problem among multiple networks that collectively predict the full solution.
For many PDEs, the operator in equation \eqref{eq:pde} can be expressed as a system of scalar equations presenting a natural separation in state space.
This case is considered in section \ref{sec:CEDNN}, where we introduce the notion of coupled EDNN (\cEDNN) for solving such systems.  
The remaining difficulties are addressed using domain decomposition, which is a well established strategy to divide the problem into a set of smaller sub-problems.
In section \ref{sec:DEDNN}, we introduce the notion of a distributed EDNN (\dEDNN):  we partition the global spatial domain into multiple elements, and solve the governing equations on each sub-domain using element-specific networks that exchange information at the element boundaries.  In section \ref{sec:correction_function}, we detail the requirements and construction of correction functions needed to communicate information between elements.
In section \ref{sec:boundary}, we describe our methodology to enforce boundary conditions.

Both \cEDNN and \dEDNN employ multiple neural networks, and are considered sub-classes of \mEDNN methods.
Figure \ref{fig:mednn} shows a schematic of these networks collectively solving a system of PDEs over many subdomains.
The domain is partitioned into multiple subdomains, each with a local vector-values solution. \cEDNN provides the methodology for the networks to coordinate their evolution within a subdomain, and \dEDNN synchronizes multiple elements by ensuring the PDE dynamics are consistent between neighboring elements.

\subsection{Coupled EDNNs}
\label{sec:CEDNN}

Consider a PDE of the form in equation \eqref{eq:pde}, where the state $\soln = \left(q_1,\dots,q_\NUMSTATE\right)$ is comprised of $\NUMSTATE \geq 2$ variables. 
There are two natural ways to approximate $\soln$ with neural networks.
One approach is to use single large neural network having a vector output, as adopted by \citet{Du_ednn} to solve the incompressible Navier-Stokes equations for example.
However, using a single network to accurately capture the entire state is challenging, in particular when the components of the solution each require different basis functions. 
To overcome this barrier, a large number of network parameters are typically needed to correctly approximate the state of the system.
As a result, the linear system \eqref{eq:linear} becomes large and computationally expensive to evolve.

An alternative and natural approach, that we term \textit{Coupled EDNN} (\cEDNN), uses multiple independent neural networks, coupled through the system of governing equations.
Each neural network is tasked with representing one of the $\NUMSTATE$ components of the state vector, so the solution ansatz \eqref{eq:ansatz} becomes,
\begin{equation}
    \hat{\soln}\left(\xdom, \weights_1(t), \dots, \weights_K(t) \right) 
    = \left(\hat{q}_1\left(\xdom,\weights_1(t)\right) ,\dots, \hat{q}_K\left(\xdom,\weights_K(t)\right) \right).
\end{equation}
In this manner, each variable within the state vector has component-specific network weights $\weights_k$, and additionally the network structure can be different.
The system of equations is then advanced according to, 
\begin{equation}
    \ddfrac{\weights_k}{\hat{q}_k} \frac{d\weights_k}{dt}
    = \mathcal{N}_k (\xdom, t, \hat{\soln}). 
\end{equation}
The parameters of each network are thus evolved separately, coupled only through their outputs which are required to evaluate the PDE operator on the right-hand side at common solution points.

\cEDNN provides several benefits:
Since networks are decoupled in parameter space, we can solve each network update separately. 
Assuming that we retain the total number of parameter $\NUMWEIGHTS$ fixed, rather than solve one large linear system \eqref{eq:optimization_linear}, we now solve $\NUMSTATE$ smaller problems one for each \cEDNN.
The estimated computational cost thus reduces from $O(\NUMX \NUMWEIGHTS^2)$ to $O(\NUMX \NUMWEIGHTS^2 / \NUMSTATE)$ \cite{bjorck}.
Gradient computations also become more efficient because they are performed using automatic differentiation on the individual networks.
Additionally, the evolution and evaluation of the individual networks can be trivially parallelized. 
Aside from the computational benefits, \cEDNN can be fine-tuned to better represent the solution. 
For example, each sub-network can utilize different architecture, or we can increase the parameterization for states where we expect to require more expressive networks.

\subsection{Distributed EDNNs}
\label{sec:DEDNN}

The idea of a \textit{Distributed EDNN} (\dEDNN) is based on a domain-decomposition strategy, where we partition the spatial domain into multiple elements and deploy an EDNN (or \cEDNN) to approximate the solution on each element.
At the element boundaries, we compute flux corrections to ensure that the dynamics are synchronized between elements and maintain a globally accurate time evolution.
Collectively, multiple EDNNs can more accurately represent complicated features of the solution compared to a single network, so we gain increased representation capability.  
Furthermore, as a finite-element methods, \dEDNN can be defined on smaller subdomains with specialized architecture for spatial regions with special characteristics such as boundary layers.  The smaller individual network sizes again leads to computational efficiency, and the method naturally lends itself to parallel implementation on distributed computing platforms.

\subsubsection*{Domain decomposition and the element sub-problem}
We focus on second-order conservation PDEs of the form,
\begin{equation}
    \label{eq:FR_pde}
    \ddt \soln + \nabla \cdot \fluxfunc\left( \soln, \nabla \soln \right) = \bm{0} 
\end{equation}
where $\fluxfunc: \Rstate\times\R^{\NUMSTATE\times\NUMDIM} \rightarrow \R^{\NUMSTATE\times\NUMDIM}$ is the flux function whose divergence is related to the PDE operator by $\operator = -\nabla\cdot\fluxfunc$.
To simplify notation, let $\nabla$ be the gradient operator in spatial variables only.
The governing equation \eqref{eq:FR_pde} can be rewritten as a system of first order equations,
\begin{equation}
    \label{eq:FR_aux_pde}
    \begin{split}
    \ddt \soln + \nabla \cdot \fluxfunc\left( \soln, \aux \right) = \bm{0},  \\
    \aux - \nabla \soln = \bm{0},
    \end{split}
\end{equation}
by introducing the auxiliary variable $\aux$.

We partition the domain $\domain$ into $\NUMELES$ elements, $\domain_e$, with boundaries $\boundary_e$.
The approximate solution on each element is $\soln_e^\dis: \domain_e \rightarrow \Rstate$, and is predicted by a \cEDNN, which immediately ensures continuity within the interior of each subdomain.
At the boundaries $\bm{\partial\Omega}_e$, however, there can be discontinuities in the approximate solutions and auxiliary variables.
We denote these discontinuous variables with superscript $\dis$.
For each element, the system of equations becomes, 
\begin{equation}
\label{eq:FR_approx_pde}
    \begin{split}
    \ddt \soln_e^\dis + \nabla \cdot \flux_e = \bm{0},  \\
    \aux_e^\dis - \nabla \soln_e = \bm{0}  & .
    \end{split}
\end{equation}
To ensure that the dynamics are consistent across element interfaces, we formulate a flux $\flux_e$ whose normal component is equal between elements.
Similarly, we formulate a solution $\soln_e$ that is continuous between elements, and thus continuous over the whole domain.

\subsubsection*{Solution correction and the auxiliary variable}
In order to obtain the continuous solution $\soln_e$, we first compute a common interface solution and a correction on each element boundary. 
This correction is then propagated to the interior of the element using a correction function.  We then compute the gradient of the corrected solution, which defines the auxiliary variable as given in \eqref{eq:FR_approx_pde}.  These steps are detailed below.  

On the shared interface between neighboring elements, we have two predictions of the solution.  
We use subscript $e-$ to denote values local to the element, and $e+$ to denote values from neighboring element(s).
The local solution $\soln_{e-}^\dis: \boundary_e \rightarrow \Rstate$ is given by the element (or \cEDNN) prediction on its own boundary $\soln_{e-}^\dis = \soln_e^\dis \vert_{\domain_e}$,
and $\soln_{e+}^\dis: \boundary_e \rightarrow \Rstate$ is from the neighbouring element (or \cEDNN).
Using a problem-specific choice of Riemann solver $\solnSolver$, we calculate the common interface solution,
\begin{equation}
    \soln_e^\star = \solnSolver(\soln_{e-}^\dis, \soln_{e+}^\dis, \bm{n}_e),
\end{equation}
where $\bm{n}_e: \boundary_e \rightarrow \Rdim$ is the outward unit normal.
Note, by design of the Riemann solver, each element common interface solution will be equal $\soln^\star_e(\xbou) = \soln^\star_{e'}(\xbou)$ at points $\xbou \in \boundary_e \cap \boundary_{e'}$ on both elements boundaries.
Thus, the necessary correction to $\soln_e^\dis$ on the boundary is the difference,
\begin{equation}
    \label{eq:FR_soln_diff}
    \soln_e^\Delta = \soln_e^\star - \soln_e^\dis \vert_{\boundary_e}.
\end{equation}
The correction to the solution on each element is then, 
\begin{equation}
    \label{eq:FR_soln_correction}
    \soln_e^\cor = \int_{\boundary_e} \soln_e^\Delta \solnCorrectFunc_e \, d \boundary_e, 
\end{equation}
which distributes the boundary correction to the interior. 
The correction function $\solnCorrectFunc_e: \boundary_e \times \domain_e \rightarrow \R$ ensures that $\soln_e^\cor$ is continuous and satisfies $\soln_e^\cor \vert_{\boundary_e} = \soln_e^\Delta$, and approximates zero in some sense (see \S\ref{sec:correction_function}). 
Therefore, the corrected solution is given by,
\begin{equation}
    \soln_e = \soln_e^\dis + \soln_e^\cor ,
\end{equation}
which is continuous across element boundaries.

The gradient of this corrected solution gives the auxiliary variable, 
\begin{equation}
    \aux_e^\dis = \nabla \soln_e^\dis + \int_{\boundary_e} \soln_e^\Delta \nabla \solnCorrectFunc_e \, d \boundary_e,
\end{equation}
where the outer product $\soln_e^\Delta \otimes \nabla \solnCorrectFunc_e$ is implied.
Note that we do not cast separate neural networks to represent the auxiliary variables.
Instead we formulate these variables as functions of the discontinuous state $\soln_e^\dis$.
This remark is particularly relevant to $\aux_{e-}^\dis$ (and $\aux_{e+}^\dis$), where we compute  
$\aux_{e-}^\dis (\xbou) = \nabla \soln_e^\dis (\xbou) + \soln_e^\Delta \nabla \solnCorrectFunc_e (\xbou,\xbou) $.
The auxiliary variable is, in general, discontinuous at element interfaces.
For second-order PDEs in general, the divergence of the flux requires the gradient of the auxiliary variable,
\begin{equation}
    \label{eq:FR_aux_grad}
    \nabla \aux_e^\dis = \nabla\nabla \soln_e^\dis + \int_{\boundary_e} \soln_e^\Delta \nabla\nabla \solnCorrectFunc_e \, d \boundary_e,
\end{equation}
where $\nabla\nabla$ denotes the Hessian of mixed, second-order partial derivatives.

\subsubsection*{Flux splitting and flux correction}
For equation \eqref{eq:FR_approx_pde}, we seek fluxes whose normal components are equal on common boundaries between elements.
Similar to the above procedure, we first compute common normal fluxes on the boundary and spread the necessary correction into the interior.
The solution on each element is then evolved according to equation \eqref{eq:FR_approx_pde}, using the divergence of the corrected flux.

In advection-diffusion problems, we separate the total flux into inviscid and viscous fluxes,
\begin{equation}
    \label{eq:flux}
    \flux_e^D = \flux_{inv,e}^\dis + \flux_{vis,e}^\dis ,
    \qquad
    \flux^\dis_{inv,e} = \fluxfunc_{inv}(\soln_e^\dis ), 
    \qquad
    \flux^\dis_{vis,e} = \fluxfunc_{vis}(\soln_e^\dis, \aux_e^\dis ).
\end{equation}
These fluxes are discontinuous at element boundaries.  
We adopt problem-specific choices of Riemann Solvers $\fluxSolver_{inv}$ and $\fluxSolver_{vis}$ to calculate the common fluxes between elements, namely the common normal-inviscid and normal-viscous fluxes,
\begin{equation}
    \flux_{inv,e}^{\perp\star} = \fluxSolver_{inv}(\soln_{e-}, \soln_{e+}, \bm{n}_e),
\end{equation}
\begin{equation}
    \flux_{vis,e}^{\perp\star}  = \fluxSolver_{vis}(\soln_{e-}, \soln_{e+}, \aux_{e-}, \aux_{e+}, \bm{n}_e).
\end{equation}
The difference between the common and discontinuous fluxes then constitutes the correction on the boundary,
\begin{equation}
    \flux_e^{\perp\Delta} = \flux_{inv,e}^{\perp\star} + \flux_{vis,e}^{\perp\star} - \flux_e^\dis \vert_{\boundary_e} \cdot \bm{n}_e .
\end{equation}
This correction is distributed on the element according to,  
\begin{equation}
    \label{eq:FR_flux_correction}
    \flux_e^\cor = \int_{\boundary_e} \flux_e^{\perp\Delta}  \fluxCorrectFunc_e \, d \boundary_e, 
\end{equation}
where $\fluxCorrectFunc_e: \boundary_e \times \domain_e \rightarrow \Rdim$ is a vector-valued correction function whose details are given in \S\ref{sec:correction_function}. 
We then formulate the corrected flux according to,
\begin{equation}
    \flux_{e} = \flux_{e}^\dis + \flux_{e}^\cor.
\end{equation}

Each element (or \cEDNN) is evolved according to the divergence of its total flux,
\begin{equation}
    \operator_e = 
    - \nabla \cdot \flux^\dis_e 
    - \int_{\boundary_e} \flux_e^{\perp\Delta}  \nabla \cdot \fluxCorrectFunc_e \, d \boundary_e.
\end{equation}
Given the splitting of the total flux, the divergence comprises two parts, 
\begin{equation}
    \nabla \cdot  \flux_{inv,e}^\dis 
    = \nabla \cdot \fluxfunc_{inv}(\soln^\dis_e , \nabla \soln^\dis_e) 
\end{equation} 
\begin{equation}
    \label{eq:divergence_vis_flux}
    \nabla \cdot \flux_{vis,e}^\dis 
    = \nabla \cdot \fluxfunc_{vis} (\soln^\dis_e, \nabla \soln^\dis_e, \aux^\dis_e, \nabla\aux^\dis_e) .
\end{equation}
The viscous flux is a function of the solution and the auxiliary variable, so its divergence, in general, requires the gradients of the solution and of the auxiliary variable. 
The computation of these divergences can be nuanced, so we provide an example from the compressible Navier-Stokes equations in \ref{sec:cns_viscous}.

In practice, we sample boundary points to approximate the solution correction \eqref{eq:FR_soln_correction} and flux correction \eqref{eq:FR_flux_correction}.
Determining the subset of points that give the most accurate and efficient evolution poses an interesting research question worthy of a separate investigation.  
In this work, we partition the domain into axis-aligned hyper-rectangular subdomains.
Then, we orthogonally project each sampling point onto all $2S$ faces to generate our set of boundary point.

\subsection{Correction functions}
\label{sec:correction_function}

Our procedure requires correction functions $\solnCorrectFunc_e$ for the solution  \eqref{eq:FR_soln_correction} and $\fluxCorrectFunc_e$ for the flux  \eqref{eq:FR_flux_correction}, to propagate the boundary data into the interior of an element. 
We first establish the requirements for one-dimensional correction functions, then introduce a new class of correction functions well-suited for \mEDNN.
We also describe a multi-dimensional extension for hyper-rectangle elements.

\subsubsection*{1D correction functions}
Consider in 1D the interval subdomain $\domain_e = \left[ x_{e,L}, x_{e,R} \right]$. 
Similar to \citet{Castonguay_13}, we define correction functions in a reference element $\domain' = \left[-1, 1\right]$ with boundary $\boundary' = \{ -1, 1 \}$.
Using the affine transformation,
\begin{equation}
    r_e(x) = \left(x - \frac{x_{e,L}+x_{e,R}}{2} \right) \frac{2}{x_{e,R}-x_{e,L}} ,
\end{equation}
physical spatial coordinates are mapped to their reference coordinates.
We formulate the reference solution correction, $\solnCorrectFunc' : \boundary' \times \domain' \rightarrow \Rdim$, as the contribution from a point on the boundary to a point within the domain.
Afterward, the solution correction function for each element is
\begin{equation}
    \solnCorrectFunc_e \left( \chi, x \right) = \solnCorrectFunc' \left( r_e(\chi), r_e(x) \right)
\end{equation} 
and has derivative
\begin{equation}
    \ddx\solnCorrectFunc_e \left( \chi, x \right) 
    = \ddfrac{x}{r_e} \ddfrac{r_e}{} \solnCorrectFunc' \left( r_e(\chi), r_e(x) \right) .
\end{equation} 
Let subscript $L$ denote a function evaluated at $r = -1$, and subscript $R$ denote a function evaluated at $r=1$.
For example, the reference correction function for the left point is $g'_L(r) = g'(-1,r)$.
In one dimension, the solution correction \eqref{eq:FR_soln_correction} simplifies to the accumulation from the left and right boundaries, 
\begin{equation}
    \soln_e^\cor(x) 
    = \soln_{e,L}^\Delta \ \solnCorrectFunc'_{L} ( r_e(x)) 
    + \soln_{e,R}^\Delta \ \solnCorrectFunc'_{R} ( r_e(x)), 
\end{equation}
where the differences $ \soln_{e,L}^\Delta, \soln_{e,R}^\Delta $ are evaluated according to equation \eqref{eq:FR_soln_diff}.

In order to maintain point-wise consistency with $\soln^C \big|_{\boundary_e} =  \soln_e^\Delta $, we require
\begin{equation}
    \label{eq:FR_correction_req}
    \begin{split}
        \solnCorrectFunc'_{L}\left(-1 \right) = 1 ,\qquad &
        \solnCorrectFunc'_{L}\left( 1 \right) = 0 ,\\
        \solnCorrectFunc'_{R}\left(-1 \right) = 1 ,\qquad &
        \solnCorrectFunc'_{R}\left( 1 \right) = 0,
    \end{split}
\end{equation}
and we expect symmetry $\solnCorrectFunc'_L(r) = \solnCorrectFunc'_R(-r)$ in order to avoid introducing bias.
For second order PDEs, we also require $\solnCorrectFunc'$ to be twice differentiable, which is required for the evaluation of the gradient of the auxiliary variable \eqref{eq:FR_aux_grad}.
We additional desire that the correction function approximates the zero function in some sense, so that the corrected solution remains a good approximation of the discontinuous solution. 
For one dimension, the flux correction function $\fluxCorrectFunc$ is scalar-valued and only needs to be once differentiable, but otherwise $\fluxCorrectFunc$ has identical requirements to $\solnCorrectFunc$.

\begin{figure} 
    \centering
    \begin{subfigure}[b]{0.30\textwidth}
        \begin{center}
        \includegraphics[height=\textwidth]{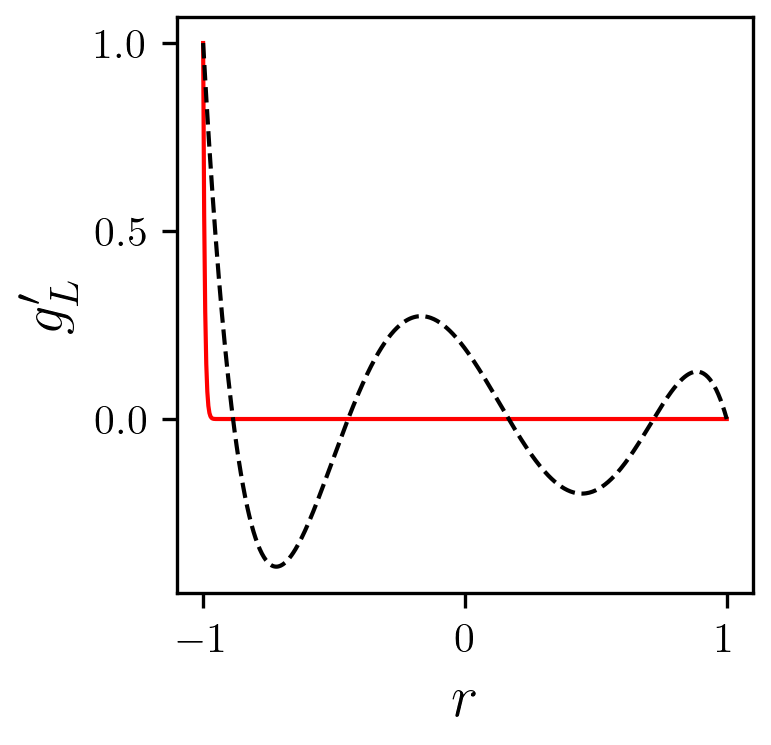}
        \end{center}
        \caption{}
        \end{subfigure}
    \hspace{0.10\textwidth}
    \begin{subfigure}[b]{0.30\textwidth}
        \begin{center}
        \includegraphics[height=\textwidth]{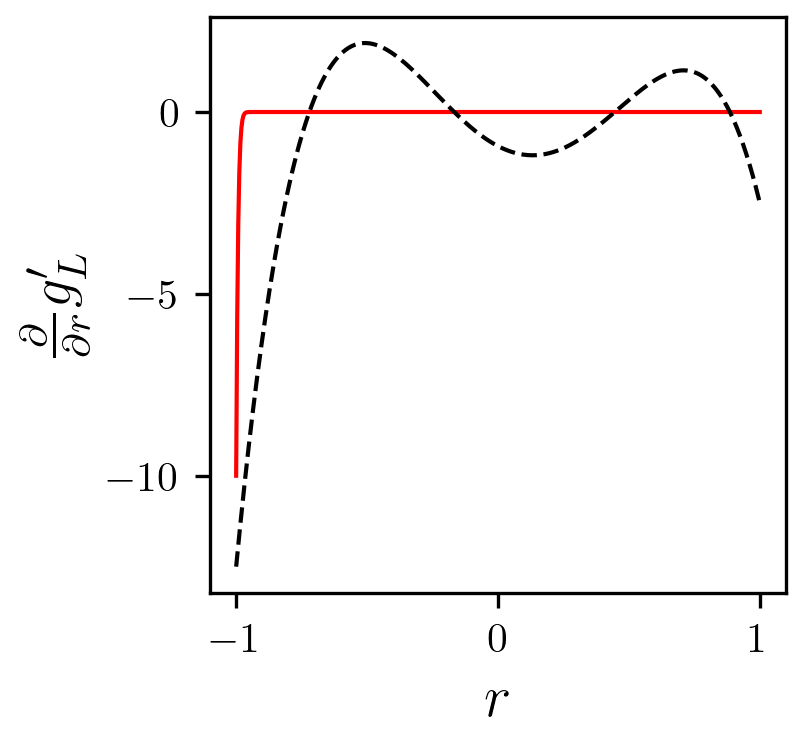}
        \end{center}
        \caption{}
    \end{subfigure}
    \caption{One-dimensional (a) correction functions and (b) their derivatives for left boundary. 
    A fifth-order polynomial $DG_5$ (black dashed) and the monomial $\mathcal{M}_{15,0.1}$ (red solid).}
    \label{fig:correction_functions}
\end{figure}

Various families of polynomial correction functions were described by \citet{Huynh_fr} and \citet{Vincent_fr} for flux reconstruction.  
% which differ in how the polynomial approximates zero.
For example, $DG_5$ refers to a 5th-order `discontinuous Galerkin' correction function, denoted $g_{DG,5}$ in \cite{Huynh_fr}, which is designed to be orthogonal to the space of polynomials of degree 3, giving some notion of being zero. 
Figure \ref{fig:correction_functions} shows $DG_5$ for the left boundary and its derivative.
In flux reconstruction, the choice of polynomial functions is natural, since the solution and flux are approximated with interpolating Lagrange polynomials.
With EDNN, however, polynomial correction functions no longer offer a good approximation to zero, since EDNN is not restricted to polynomial function spaces. 
Also, the EDNN sampling points may be anywhere within the element and can be changed dynamically during the evolution, which renders a polynomial correction function undesirable especially since its derivative is non-zero far from the boundary (for $DG_5$, figure \ref{fig:correction_functions} shows $\ddr\solnCorrectFunc'_{e,L}$ has full support and is non-zero near $r=1$).  
Intuitively, the correction from one boundary should not have a sizeable effect on points far from that boundary.  This issue is particularly relevant to \dEDNN, where the element sizes can be large when the networks have high expressivity and can capture the solution over large portions of the global domain.   

\subsubsection*{Monomial correction functions}
For use with EDNN, we construct correction functions that better approximate zero in more general function spaces. 
With order parameter $\mathcal{P}$ and width parameter $\mathcal{W}$, we devise the monomial class $\mathcal{M}_{\mathcal{P},\mathcal{W}}$ of correction functions,
\begin{equation}
\solnCorrectFunc'_{L}(r) =  
    \begin{cases}
    \left( \frac{r+1-\mathcal{W}}{-\mathcal{W}}\right)^\mathcal{P} & \text{if } r \leq -1+\mathcal{W} \\
    0 & \text{otherwise } 
    \end{cases} 
, \qquad
\solnCorrectFunc'_{R}(r) =  
    \begin{cases}
    \left( \frac{r-1+\mathcal{W}}{\mathcal{W}}\right)^\mathcal{P} & \text{if } r \geq 1-\mathcal{W} \\
    0 & \text{otherwise } 
    \end{cases} . 
\end{equation}
These functions are a $\mathcal{P}^{\textrm{th}}$ order monomial in the interval of width $\mathcal{W}$ adjacent to the respective boundary point, and are zero everywhere else.
In figure \ref{fig:correction_functions}, we show that monomial $\mathcal{M}_{5,0.1}$ and its derivative quickly decay to zero.  
Therefore, the monomials satisfy our desire to only impact the interior near the boundary.
In the limit as $\mathcal{P} \rightarrow \infty$ or $\mathcal{W} \rightarrow 0$, this correction function approaches an indicator function which is zero in the interior of the domain.
Our preliminary testing showed that, for EDNN, these monomials yield more accurate predictions than polynomial correction functions.

\subsubsection*{Extension to higher dimensions}
To construct correction functions in $\NUMDIM > 1$ dimensions, we consider partitions of axis-aligned hyper-rectangles. 
We then cast a series of one-dimensional sub-problems, as detailed in the work by \citet{Huynh_fr} and also \citet{Castonguay_13}.
Formally, the solution correction function becomes,
\begin{equation}
    \solnCorrectFunc_e (\xbou, \xdom) = 
    \begin{cases}
        g_{L}'(r_{e,s}(\xdom)) & 
        \text{if } \exists~s : \xbou = \xbou_{e,s,L} \\
        g_{R}'(r_{e,s}(\xdom)) & 
        \text{if } \exists~s : \xbou = \xbou_{e,s,R} \\
        0         & \text{otherwise }
    \end{cases}, 
\end{equation}
where subscript $s$ denotes a dimension in the spatial coordinates.
In the above expression, $\xbou_{e,s,L}$ and $\xbou_{e,s,R}$ are the orthogonal projections of $\xdom$ onto the left and right faces of the $s^{\textrm{th}}$ dimension.
The solution correction reduces to the contributions from the $2S$ boundary points,
\begin{equation}
    \soln_e^\cor(\xdom) = \sum_{s=1}^\NUMDIM
      \soln_e^\Delta(\xbou_{e,s,L}) \solnCorrectFunc_L'\left(r_{e,s}(\xdom) \right) 
    + \soln_e^\Delta(\xbou_{e,s,R}) \solnCorrectFunc_R'\left(r_{e,s}(\xdom) \right) .
\end{equation}
The flux correction function is extended similarly, 
\begin{equation}
    \fluxCorrectFunc (\xbou, \xdom) = 
    \begin{cases}
        h_{L}'(r_{e,s}(\xdom)) \bm{e}_s & 
        \text{if } \exists~s : \xbou = \xbou_{e,s,L} \\
        h_{R}'(r_{e,s}(\xdom)) \bm{e}_s & 
        \text{if } \exists~s : \xbou = \xbou_{e,s,R} \\
        0         & \text{otherwise }
    \end{cases},
\end{equation}
where $h_{L}'$ and $h_{R}'$ are the one-dimensional flux correction functions and $\bm{e}_s$ is the  unit basis vector in the $s^{\textrm{th}}$ direction.
The flux correction is therefore, 
\begin{equation}
    \flux_e^\cor(\xdom) = \sum_{s=1}^\NUMDIM
      \flux_e^{\perp\Delta}(\xbou_{e,s,L}) h_{L}'(r_{e,s}(\xdom)) \bm{e}_s 
    + \flux_e^{\perp\Delta}(\xbou_{e,s,R}) h_{R}'(r_{e,s}(\xdom)) \bm{e}_s .
\end{equation}

\subsection{Boundary conditions}
\label{sec:boundary}

\subsubsection*{Periodic boundary conditions}
In \mEDNN, the flux correction procedure can be adopted to enforce periodic boundary conditions in a straightforward fashion. 
Element interfaces on the boundary $\boundary$ are treated similarly to interfaces within the interior of the domain.
The same approach can also be used when a single EDNN (or \cEDNN) is adopted through the entire solution domain.  

In some of our numerical experiments, we wish to compare \mEDNN to a reference single-element solution that is independent of the boundary-correction algorithm.
In these instances, the single-element solution will adopt the feature expansion approach described by \citet{Yazdani}, which maps hyper-rectangular domains to the torus $\mathbb{T}^\NUMDIM$.
This feature layer thus restricts the single-element EDNN prediction to the space of periodic functions.

\subsubsection*{Time-dependent Dirichlet boundary conditions}
Consider $\NUMBOUND$ non-intersecting objects with known and possibly time-dependent geometry $\domain_{b}$ and state $\soln_{b}$. 
We formulate a modified operator,
\begin{equation}
    \label{eq:solid_rhs}
    \ddfrac{t}{\hat{\soln}} = 
    \begin{cases}
        \bcWeightSolid \left(\soln_{b} - \hat{\soln} \right)  & \text{if } \exists~b : \xdom \in \domain_{b} \\
        \operator(\hat{\xdom, t, \soln})                                & \text{otherwise} \\
    \end{cases},
\end{equation}
where $\bcWeightSolid$ is a parameter that determines the rate of relaxation. 
For every sampling point $\xdom_i \in \domain$, we first determine whether it intersects an object geometry, in which case we drive EDNN to match the object state.
There is no guarantee, however, that the sampling points from $\domain$ intersect with every object.
For example, in Couette flow (\S\ref{sec:couette}), the top and bottom walls are a measure zero subset of the domain. 
To ensure that EDNN tracks the objects, we can sample additional points directly from $\domain_{b}$.
In a \mEDNN context, we then determine which element subdomain contains these points, and include them to the associated EDNN linear system \eqref{eq:linear}.

This method only weakly enforces boundary conditions. 
EDNN will be driven to the object state at a rate of ${1}/{\bcWeightSolid}$.
Another interpretation is that the modified operator provides a source term to the governing equation. 
This methodology can also be used for stationary and moving solid objects immersed in the fluid.
With a slight modification, 
$\partial \hat{\soln} / \partial t = \operator(\hat{\soln}) + c_s(\xdom) \left(\soln_{b} - \hat{\soln} \right)$, 
we can implement sponge zones \citep{Bodony} to either absorb incident disturbances or to introduce them into the domain.

\section{Numerical experiments}
\label{sec:experiments}

In this section, we demonstrate \mEDNN capability and accuracy by solving canonical PDE problems that are relevant to fluid dynamics, including advection, diffusion, and the Navier-Stokese equations.  
To demonstrate \cEDNN, where multiple netowrks cooperatively solve coupled equations on a single element, we directly solve the compressible Navier-Stokes equations for the Couette-flow problem.  Starting from an initial condition far from the steady state, we show that \cEDNN converges to an analytical reference solution.
To examine the performance of \dEDNN over spatially distributed elements, we first solve the one-dimensional linear advection equation, which exercises the inviscid fluxes, and show accurate prediction over many subdomains. 
We then exercise the viscous fluxes across \dEDNN elements by solving the two-dimensional time-dependent heat equation.
Lastly, we utilize the full \mEDNN framework, leveraging both state-space (\cEDNN) and spatial decomposition (\dEDNN) to solve the Navier-Stokes equations for Taylor-Green vortices.

Unless otherwise stated, every neural network utilizes input and output scaling layers, as described by \citet{Yazdani}.
The input scaling achieves the mapping $\domain_e \rightarrow [-1,1]^\NUMDIM$, mimicking mapping to a reference, or standard, element.
The output scaling layer maps $[-1,1] \rightarrow [q_{k,L}, q_{k,R}]$ where $q_{k,L}$ and $q_{k,R}$ are the minimum and maximum of the initial state.
For our tests, we use the hyperbolic tangent (tanh) activation function and more than one hidden layers.
We recognize the idea by \citet{Anderson_22} and \citet{Bruna_23} that high accuracy can be achieved by judiciously selecting a network architecture such that the solution ansatz matches the solution of the PDE, in a single-layer shallow network.
Ultimately, however, we wish to utilize the expressivity of deep neural networks to solve complex PDE problems, and for this purpose we adopt general network architectures and without exploiting \emph{a priori} knowledge of the solution form.
For this same reason, in our \dEDNN configurations, we choose relatively large subdomains, or elements, compared to flux reconstruction experiments.  For example, in the one-dimensional linear advection case, we choose subdomains that can fully support the full waveform, while flux reconstruction employs a large number of small elements so that the solution is approximately polynomial within an element.

\subsection{Couette flow}
\label{sec:couette}

We start by demonstrating \cEDNN for solving the compressible Navier-Stokes equations.   We denote dimensional variables with superscript $\star$, and non-dimensional quantities are scaled by 
the characteristic length scale $L^\star$,
the speed of sound $c_\infty^\star$,
a reference density $\rho_\infty^\star$, 
and the specific heat capacity at constant pressure $C_{p,\infty}^\star$.
This choice yields the Reynolds number $ Re = {\rho_\infty^\star c_\infty^\star L^\star}/{\mu_\infty^\star}$ and Prandtl number $Pr = {C_{p,\infty}^\star \mu^\star}/{\kappa^\star}$, 
where $\mu^\star$ is the dimensional viscosity and $\kappa^\star$ is the dimensional thermal conductivity.
The governing PDE in non-dimensional, conservative form becomes,
\begin{equation}
    \ddfrac{t}{\soln} + {\nabla\cdot} \fluxfunc_{inv} + {\nabla\cdot} \fluxfunc_{vis} = 0
\end{equation}
where $\soln = \left[\rho, \rhou_s, \rhoe \right]^\top$ are the conservative variables: density $\rho$, momenta $\rhou_s$ in spatial directions $s = \{1, 2, 3$\}, and total energy $\rhoe$.
The inviscid and viscous fluxes in direction $i$ are
\begin{equation}
    \label{eq:cns_flux}
    \fluxfunc_{inv,i}
    =
    \begin{bmatrix}
    \rho u_i  \\
    \rhou_s u_i + p\delta_{si} \\
    u_i \left( \rhoe + p \right) \\
    \end{bmatrix}
    , \qquad
    \fluxfunc_{vis,i}
    =
    \begin{bmatrix}
    0 \\
    - \tau_{si} \\
    - u_k \tau_{ki} + \theta_i  \\
    \end{bmatrix},
\end{equation}
with summation over repeated indices, and $\delta$ denotes the Kronecker delta.
The primitive velocities, pressure, temperature, and first viscosity are,
\begin{equation}
    u_i = \frac{\rhou_i}{\rho}
    ,\qquad
    p = \left( \gamma-1 \right) \left( \rhoe - \frac{1}{2} \rhou_k \, u_k\right) 
    ,\qquad
    T = \frac{p}{R_g \rho}
    ,\qquad
    \mu = \left( \left( \gamma -1 \right) T \right)^\alpha,
\end{equation}
where $\gamma$ is the ratio of specific heats, $R_g$ is the ideal gas law constant, and $\alpha$ is the exponent of the viscosity power law. 
The viscous stress tensor and the heat flux tensor are,
\begin{equation}
    \tau_{ki}
    = \frac{\mu}{Re} \left( \ddfrac{x_k}{u_i} + \ddfrac{x_i}{u_k} \right) 
    + \frac{1}{Re}\left(\mu_b - \frac{S-1}{S }\mu \right) \ddfrac{x_l}{u_l} \delta_{ki}
    , \qquad
    \theta_i = - \frac{\mu}{Re Pr} \ddfrac{x_i}{T}
\end{equation}
where $\mu_b$ is the bulk viscosity.
In our numerical experiments, we take $\gamma = 1.4$, $\mu_b = 0.6$, and $Pr = 0.72$. 
The remaining problem-specific parameters for the viscosity power law $\alpha$ and the Reynolds number $Re$ are defined for each problem. 

\begin{figure}
     \centering
     \begin{subfigure}[b]{0.30\textwidth}
         \centering
         \includegraphics[height=\textwidth]{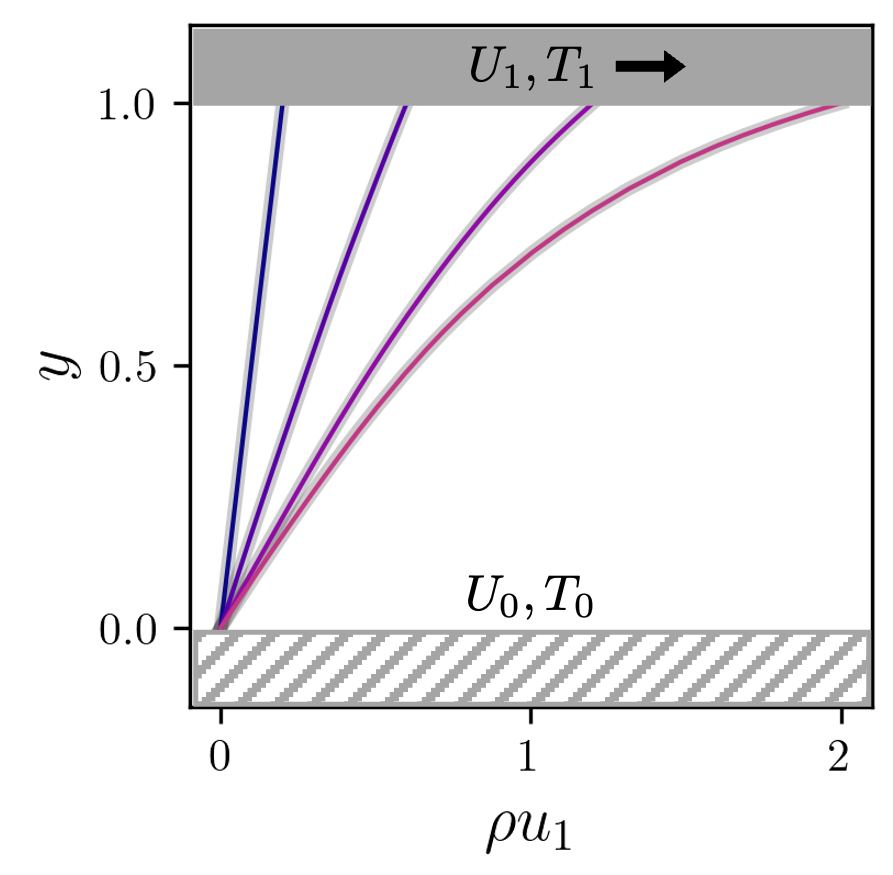}
         \caption{}
         \label{fig:couette_result_a}
     \end{subfigure}
     \begin{subfigure}[b]{0.30\textwidth}
         \centering
         \includegraphics[height=\textwidth]{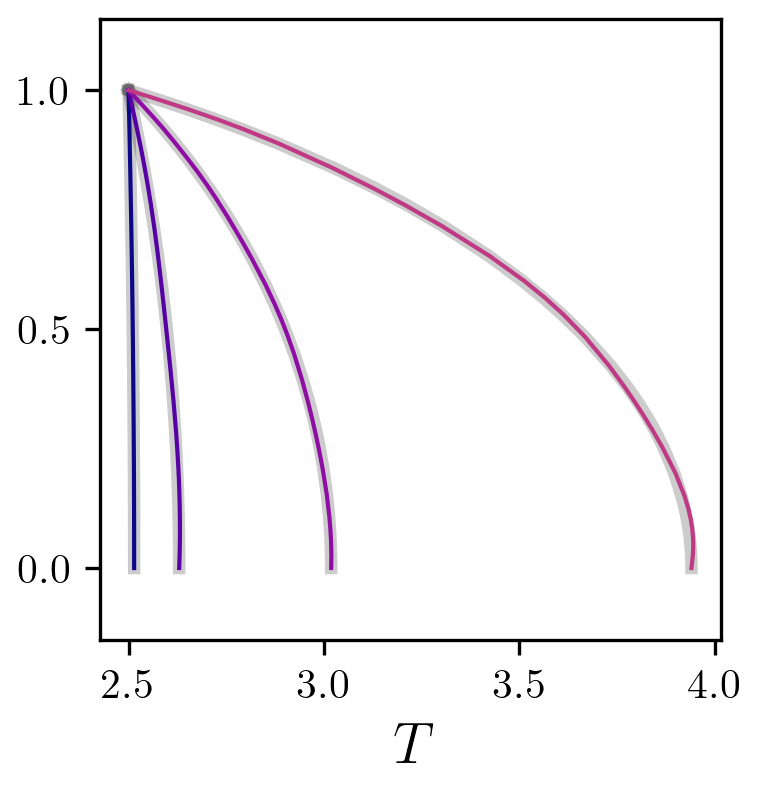}
         \caption{}
         \label{fig:couette_result_b}
     \end{subfigure}
     \begin{subfigure}[b]{0.30\textwidth}
         \centering
         \includegraphics[height=\textwidth]{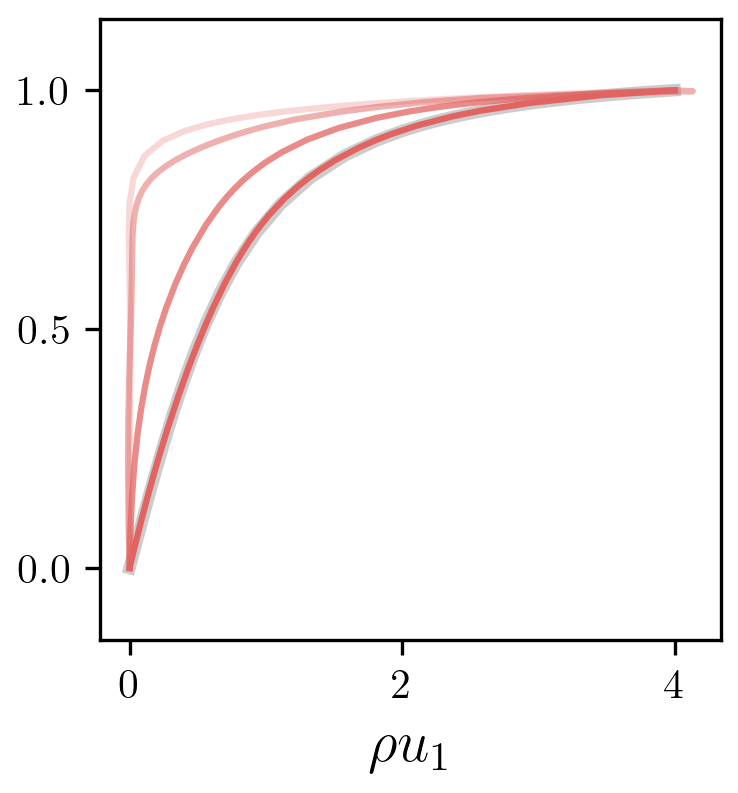}
         \caption{}
         \label{fig:couette_result_c}
     \end{subfigure}
    \caption{\cEDNN solution of Couette flow. 
    Analytical (grey) and \cEDNN prediction (colored) of steady-state (a) $\rhou_1$ and (b) temperature for Mach numbers $Ma = \{0.2, 0.5, 1.2, 2.0\}$. 
    (c) Analytical steady-state (grey) and \cEDNN prediction (colored) of $\rhou_1$ for $Ma = 4$, $Re = 100$ at times $t= \{0, 0.1, 1, 10\}$.}
    \label{fig:couette_result}
\end{figure}

Using these governing equations, we solve the classic Couette flow problem in a two-dimensional domain $\Omega = [0,1]^2$.
The top plate is a moving no-slip isothermal wall with streamwise velocity $U_1$ and temperature $T_1$, and the bottom plate is a stationary isothermal no-slip wall with velocity $U_0$ and temperature $T_0$.
Figure \ref{fig:couette_result_a} shows a schematic of the flow configuration. 

For this problem, we treat the entire domain as a single element, and adopt \cEDNN for the solution of the compressible Navier-Stokes equations.  We deploy four networks, one for the continuity equation, two for the momentum equations in the streamwise and wall-normal directions, and one for the energy equations.  We start from initial conditions that are far from the steady state, and our goal is verify whether \cEDNN converges to the steady-state analytic solution.

With viscosity power law exponent $\alpha = 1$, the Navier-Stokes equations permit the implicit analytical solution,
\begin{align}
%    \begin{cases}
        \frac{u_1}{U_1} + \frac{1}{2} Pr (\gamma -1) U_1^2 \left(\frac{u_1}{U_1} - \frac{1}{3}\left(\frac{u_1}{U_1} \right)^3 \right) &= \left(1 + \frac{1}{3} Pr (\gamma - 1) U_1^2\right)y \nonumber \\
        \frac{T}{T_1} = 1  + \frac{1}{2} Pr (\gamma - 1) \left(U_1^2 - u_1^2 \right), \qquad
        p &= R_g T_1, \qquad
        u_2 = 0
 %   \end{cases}
\end{align}
from which we can derive all other fluid variables, in particular the conserved variables, such as $\rhou_1$ shown in figure \ref{fig:couette_result}.
Choosing a top wall temperature $T_1 = \frac{1}{\gamma - 1}$, we will consider multiple cases each with different value of the top-wall speed, $U_1= \{0.2, 0.6, 1.2, 2.0, 4.0\}$.
Note that the velocity is normalized by the speed of sound, and therefore the values of the top-wall speeds are also the relevant Mach numbers of the problem.
The bottom wall boundary conditions become $U_0=0$ and $T_0 = T_1  + \frac{1}{2} Pr U_1^2$.
The Reynolds number is fixed to $Re=100$, which is again based on the speed of sound and hence spans from $Re_U\left(\equiv {\rho_\infty^\star U_1^\star L^\star}/{\mu_\infty^\star}\right) = 20$ to $Re_U=400$ based on the speed of the top wall.  
The initial condition starts from a large deviation away from the steady state, specifically we adopt the initial streamwise velocity profile is $u_1 = U_1 y^{20}$. 

Since the problem is constant in $x$, we use the input feature expansion layer $(x,y) \mapsto (y)$ then map $y$ using the standard scaling layer.
Each state neural network uses two layers and twenty neurons per layer, and the \cEDNN is evaluated at $32$ points adaptively sampled using the technique in \ref{sec:adapt_sampling}.
We march with RK4 and $\Delta t=10^{-6}$.
The top and bottom wall boundary conditions were enforced with 
$\bcWeightSolid = 1/\Delta t$.

In all cases, the \cEDNN prediction converged to the analytical steady-state solution which is shown using the grey lines in figure \ref{fig:couette_result}.
We show the converged $\rhou_1$ in figure \ref{fig:couette_result_a} and temperature in figure \ref{fig:couette_result_b}.
In figure \ref{fig:couette_result_c}, we show the time evolution of momentum $\rhou_1$ as it approaches and converges to the steady state, for the Mach 4 case. 
For this condition, the L2 error normalized by the norm of the steady-state solution at $t=20$ is less than $O(10^{-3})$ 
for all variables.

\subsection{Linear advection}
\label{sec:advection}

For our first test of a distributed set of networks over space, or \dEDNN, we simulate the 1D linear advection of a narrow Gaussian over an interval $\domain = [0,\ell ]$, where the wave is generated by the condition at at the left boundary.  The governing equation and initial and boundary conditions are, 
\begin{equation}
    \ddfrac{t}{q} + \nabla \cdot \left( u q \right) = 0 ,\qquad
    q(x,0) = q_0(x) ,\qquad
    q(0,t) = q_0(-ut) ,\qquad
    q_0(x) = \exp\left( \frac{-(x-x_0)^2}{2\Sigma} \right),
\end{equation}
which have the analytical solution, 
\begin{equation}
    q(x,t) = q_0(x-ut).
\end{equation}
For our numerical experiment, the wavespeed is constant $u=1$, and the center of the initial Guassian is $x_0 = -0.5$ and its variance is $\Sigma = 0.01$, and thus the initial Gaussian is outside the computational domain.  
The domain is partitioned into unit subdomains $\domain_e = [e-1, e]$, where $e=1\dots \ell$.  Each the governing equation is solved on each the of the $\ell$ sub-domain using a separate EDNN which has four hidden layers and twenty neurons per layer.
Each EDNN is evolved with 250 points distributed uniformly in space.
We march with RK4 and $\Delta t=10^{-6}$.
The inflow boundary conditions was enforced with $1$ point and $\bcWeightSolid = 1/\Delta t$.

\begin{figure} 
    \centering
    \begin{subfigure}[b]{0.245\textwidth}
        \begin{center}
        \includegraphics[height=1.05\textwidth]{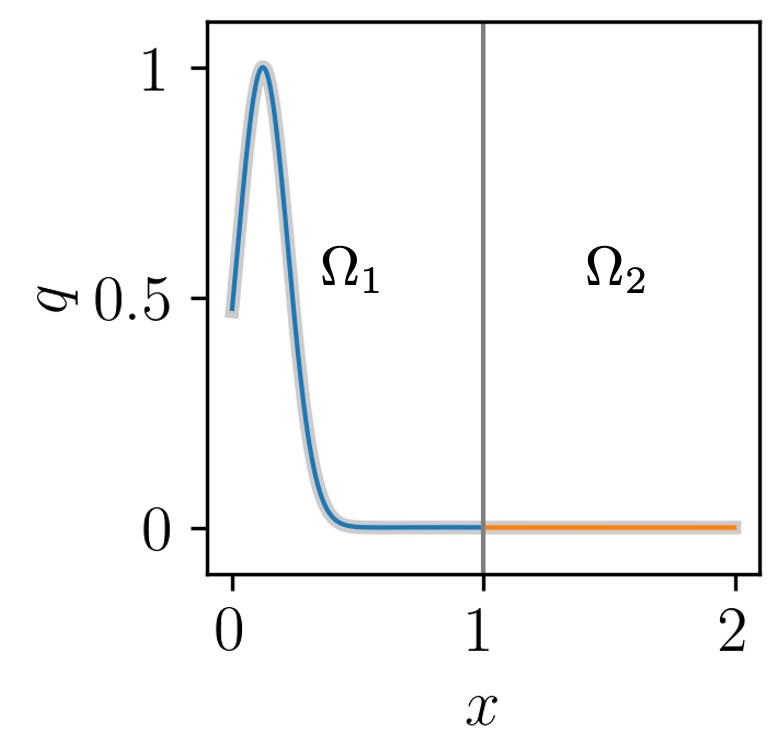}
        \end{center}
        \caption{}
    \end{subfigure}
    \hfill
    \begin{subfigure}[b]{0.245\textwidth}
        \begin{center}
        \includegraphics[height=1.05\textwidth]{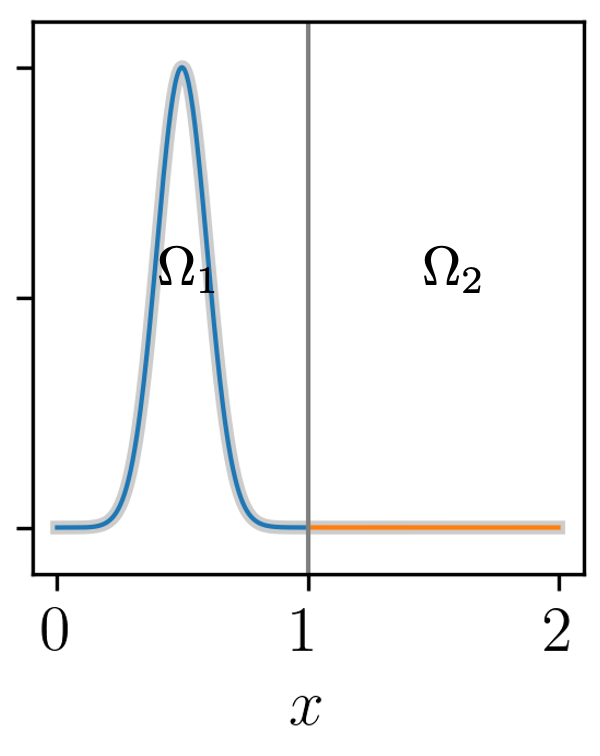}
        \end{center}
        \caption{}
    \end{subfigure}
    \hfill
    \begin{subfigure}[b]{0.245\textwidth}
        \begin{center}
        \includegraphics[height=1.05\textwidth]{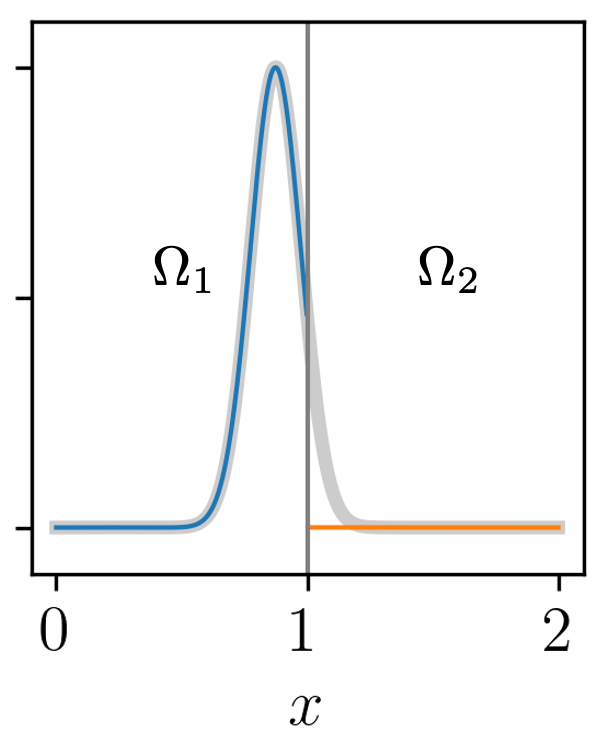}
        \end{center}
        \caption{}
    \end{subfigure}
    \hfill
    \begin{subfigure}[b]{0.245\textwidth}
        \begin{center}
        \includegraphics[height=1.05\textwidth]{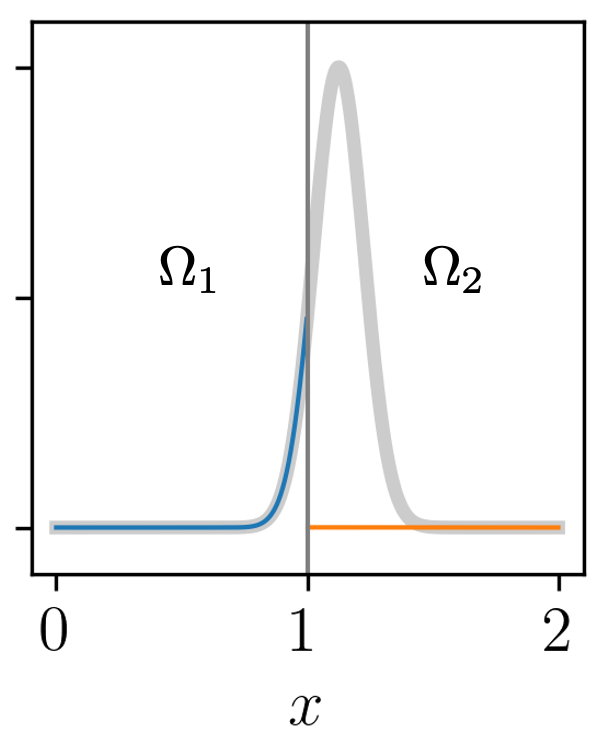}
        \end{center}
        \caption{}
    \end{subfigure}
    \caption{Predictions of 1D advection equation for $\ell=2$ without correction. Analytical solution (grey) and \dEDNN prediction (colored) at time (a) $t=0.625$, (b) $t=1$, (b) $t=1.375$, and (d) $t=1.625$ }
    \label{fig:adv_results_L02}
\end{figure}

We demonstrate the necessity of the \dEDNN procedure for the treatment of the inter-element interfaces using a configuration that fails.
We take $\ell=2$ and neglect the flux correction ($\flux_e^\cor=0$), so the total flux $(\flux_e = \flux_e^\dis$) uses only local information.
Figure \ref{fig:adv_results_L02} shows EDNN prediction for this first configuration.
Our boundary condition method successfully enforces the time-dependent inflow condition, and we see the wave correctly materialize and advect within element 1.
The wave, however, fails to cross over to the second element, $\domain_2$.

\begin{figure} [ht]
    \centering
    \begin{subfigure}[b]{0.23\textwidth}
        \begin{center}
        \includegraphics[height=1.0\textwidth]{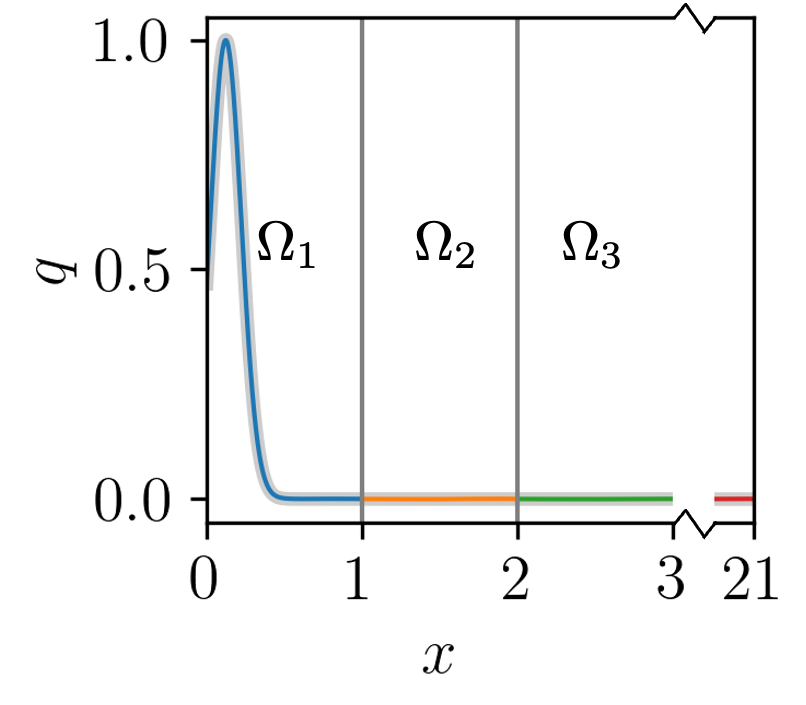}
        \end{center}
        \caption{}
    \end{subfigure}
    \hspace{0.05\textwidth}
    \begin{subfigure}[b]{0.23\textwidth}
        \begin{center}
        \includegraphics[height=1.0\textwidth]{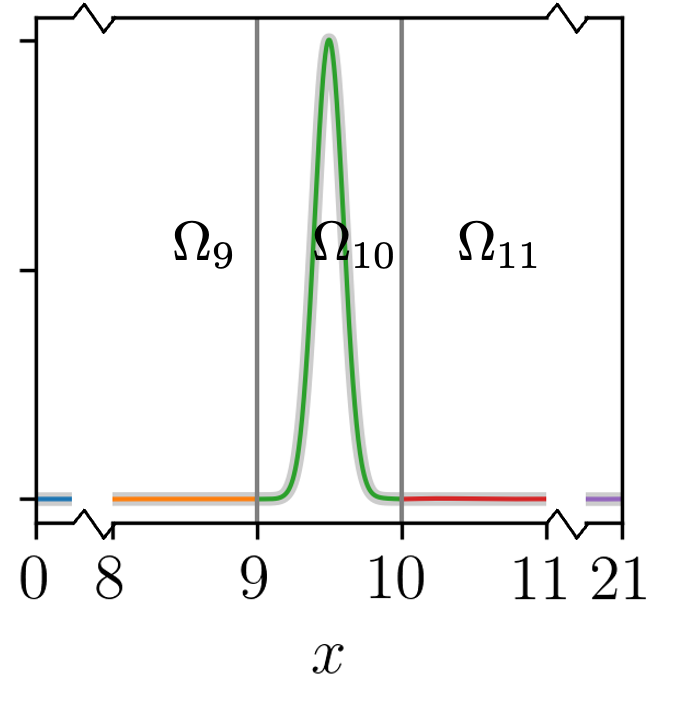}
        \end{center}
        \caption{}
    \end{subfigure}
    \begin{subfigure}[b]{0.23\textwidth}
        \begin{center}
        \includegraphics[height=1.0\textwidth]{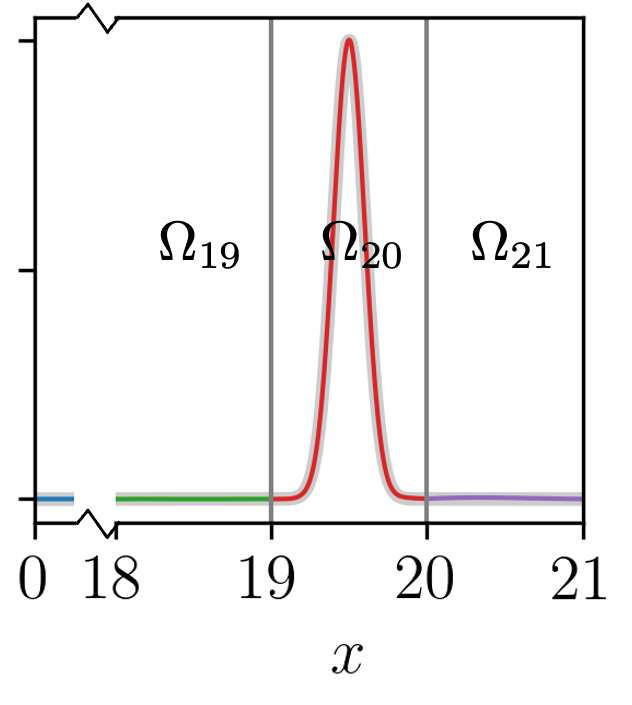}
        \end{center}
        \caption{}
    \end{subfigure}
    \begin{subfigure}[b]{0.23\textwidth}
        \begin{center}
        \includegraphics[height=1.0\textwidth]{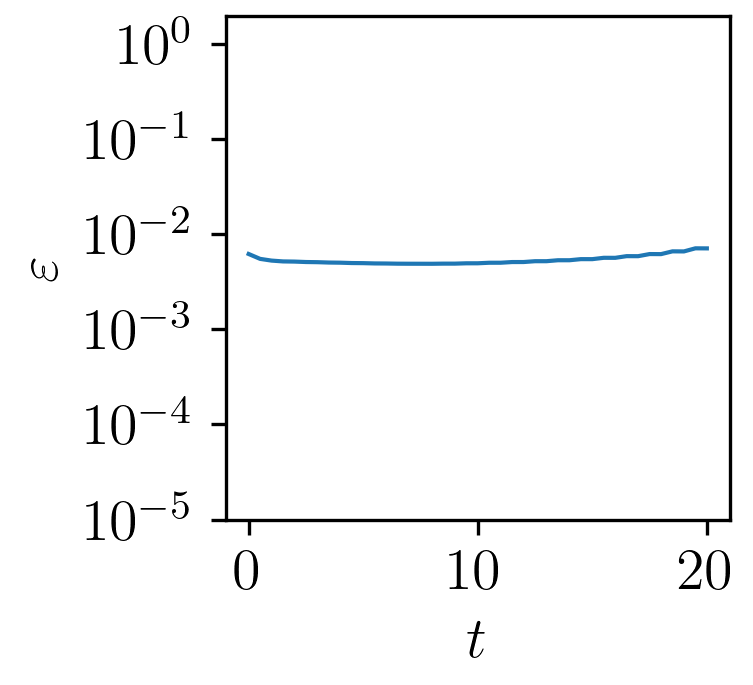}
        \end{center}
        \caption{}
        \label{fig:adv_results_L22_error}
    \end{subfigure}
    \caption{Predictions of 1D advection equation for $\ell=21$ with correction. Analytical solution (grey) and \dEDNN prediction (colored) at time (a) $t=0.625$, (b) $t=10$, and (c) $t=20$.
    (d) Instantaneous error $\varepsilon$ over time.
    }
    \label{fig:adv_results_L22}
\end{figure}

We proceed to demonstrate successful advection by adopting the interface-correction procedure.
We consider $\ell=21$ elements, adopt an upwind common inviscid flux (\ref{sec:riemann_solver}), and a monomial flux correction function $\mathcal{M}_{15,0.1}$.
Figure \ref{fig:adv_results_L22} shows that the second configuration successfully tracks the wave across element boundaries. 
To provide a quantitative measure of accuracy, we track the instantaneous error normalized by the norm of the Gaussian, 
\begin{equation}
    \varepsilon(t) = \frac{\Vert \hat{q}(x,t)-q(x,t)\Vert_2}{\Vert q_0(x) \Vert_2}.
\end{equation}
Figure \ref{fig:adv_results_L22_error} shows that the error remains commensurate with the initial training error, and the final error at $t=20$ is $\varepsilon = 7.11\times10^{-3}$.

\subsection{Linear diffusion}
\label{sec:diffusion}

For our next test, we solve the two-dimensional linear diffusion equation on a periodic domain $\domain = \left[-\pi, \pi \right]^2$, starting from a sinusoidal initial condition. The governing equations, and boundary and initial conditions are, 
\begin{equation}
    \ddfrac{t}{q} + \nabla \cdot \left( -\nu \nabla q \right) = 0,\quad
    q(-\pi,y,t) = q(\pi,y, t) ,\quad
    q(x,-\pi,t) = q(x,\pi, t) ,\quad
    q(x,y,0) = \sin x \sin y 
\end{equation}
where $\nu=1$ is the diffusivity coefficient. 
The analytical solution is given by
\begin{equation}
    q(x,y,t) = e^{-2\nu t} \sin x \sin y.
\end{equation}

\begin{figure} 
    \centering
    \begin{subfigure}[b]{0.23\textwidth}
        \begin{center}
        \includegraphics[height=\textwidth]{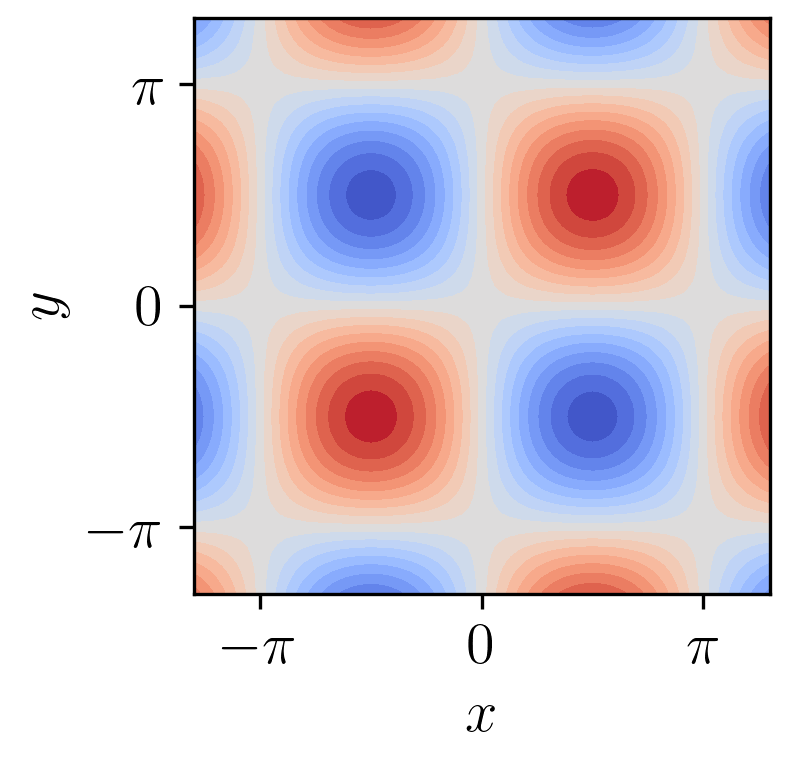}
        \end{center}
        \caption{}
    \end{subfigure}
    \begin{subfigure}[b]{0.23\textwidth}
        \begin{center}
        \includegraphics[height=\textwidth]{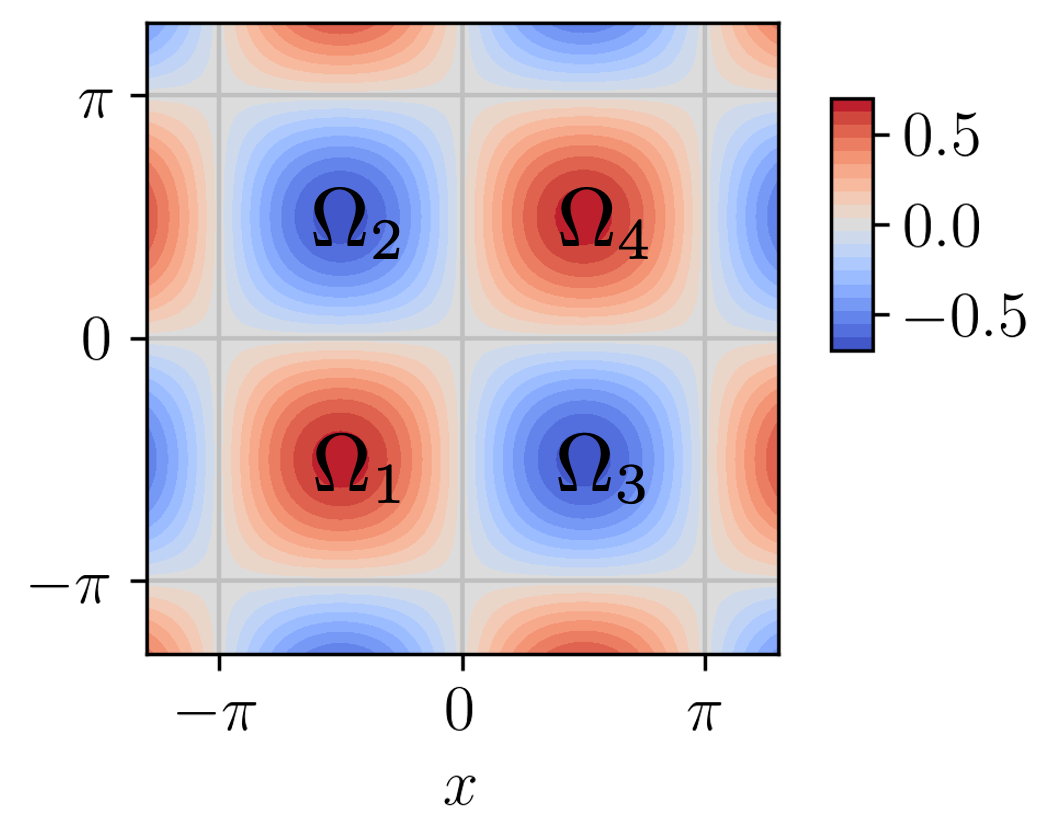}
        \end{center}
        \caption{}
        \label{fig:dif_result_ednn}
    \end{subfigure}
    \hspace{0.05\textwidth}
    \begin{subfigure}[b]{0.23\textwidth}
        \begin{center}
        \includegraphics[height=\textwidth]{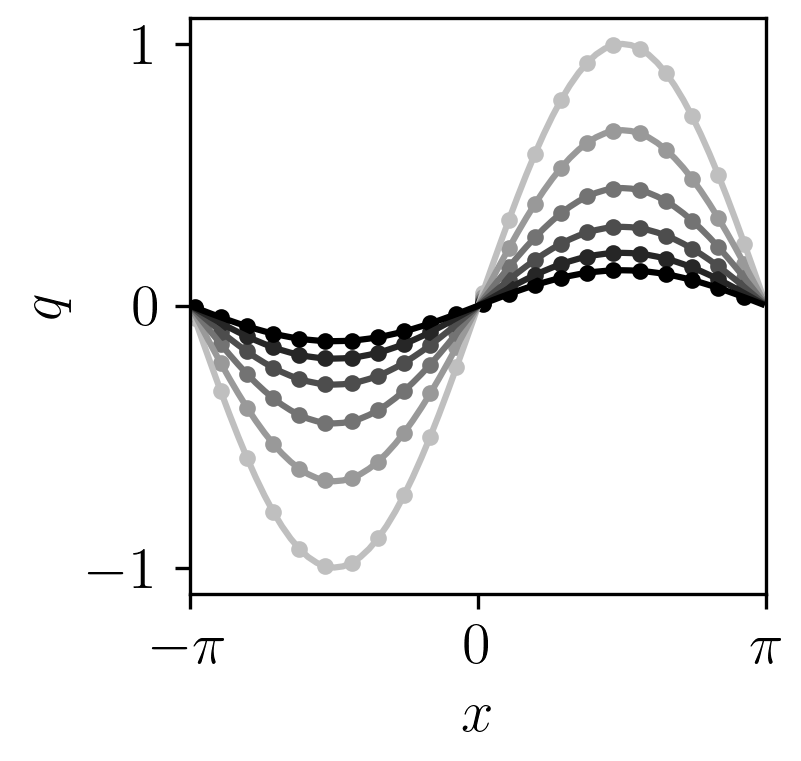}
        \end{center}
        \caption{}
        \label{fig:dif_result_compare}
    \end{subfigure}
    \begin{subfigure}[b]{0.23\textwidth}
        \begin{center}
        \includegraphics[height=\textwidth]{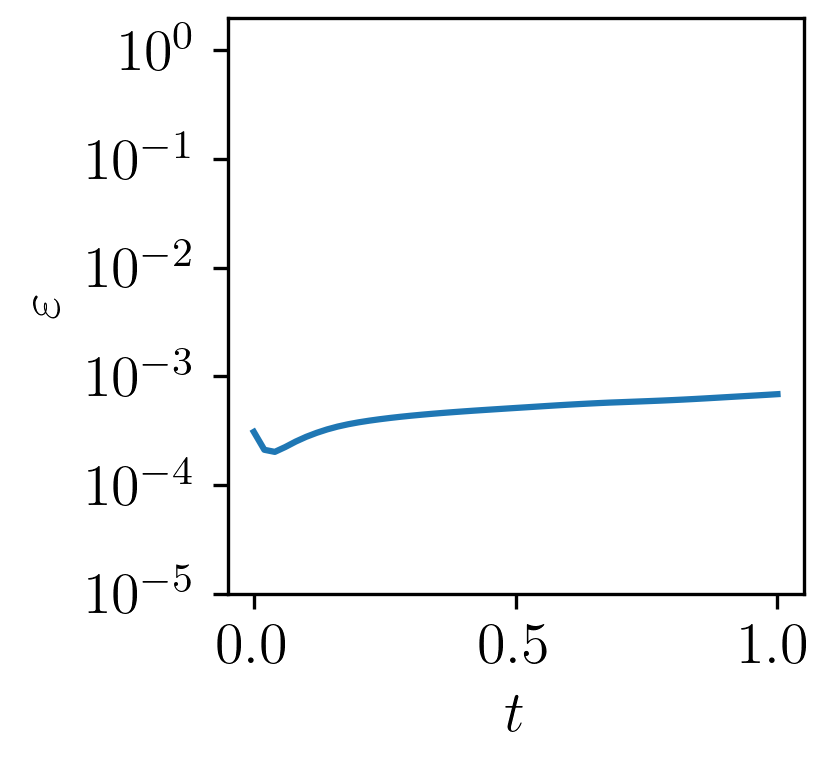}
        \end{center}
        \caption{}
        \label{fig:dif_result_error}
    \end{subfigure}
    \caption{Solution for 2D diffusion equation. 
            Contours of (a) analytical and (b) \mEDNN solution at $t=0.2$ and grey lines marking element interfaces.
            (c) Comparison between analytical (dots) and \mEDNN solution (lines) at $y=\pi/2$ and times $t= \{ 0.0,0.2,0.4,0.6,0.8,1.0 \}$.
            (d) Instantaneous error $\varepsilon$ over time.}
    \label{fig:dif_result}
\end{figure}

We adopt a spatially distributed set of networks, or \dEDNN, and partition the domain into $2\times 2$ equal elements as shown in figure \ref{fig:dif_result_ednn}. 
These partitions were selected so that element interfaces lie on the regions with maximum heat flux.
The interfaces within the domain interior and on the global periodic boundaries exchange the solution and flux corrections.
We evaluate the common solution and common normal viscous fluxes using the LDG method (\ref{sec:riemann_solver}), and both the solution and flux correction functions are the monomial $\mathcal{M}_{3,0.25}$.
Each of the four EDNNs is comprised of $4$ hidden layers and $10$ neurons per layer, and each sub-domain is sampled at $33\times33$ points distributed uniformly.
We march with RK4 and a stepsize $\Delta t = 0.01$. 

The \dEDNN predictions are compared to the analytical solution in figure \ref{fig:dif_result}, and shows very good agreement.
For visualization only, we extend the domain to show the periodicity of the \dEDNN prediction.
Figure \ref{fig:dif_result_compare} shows a comparison along the horizontal line $y=\pi/2$ at various times.
Globally, the collective prediction by the distributed \dEDNN networks matches the reference solution.
Since the temperature decays as a function of time, we evaluate the instantaneous error, 
\begin{equation}
    \label{eq:error}
    \varepsilon(t) = \frac{\Vert \hat{q}(x,t)-q(x,t)\Vert_2}{\Vert q(x,t) \Vert_2},
\end{equation}
which is reported in figure \ref{fig:dif_result_error}.  The error remains commensurate with its value from the representation of the initial condition, and at the final time, $\varepsilon = 6.86\times10^{-4}$.

\subsection{Taylor-Green vortices}

To demonstrate the full, coupled and distributed \mEDNN framework, with multiple states and multiple elements, we simulate two-dimensional Taylor-Green vortices. 
The system of PDEs is the compressible Navier-Stokes equations, as described in section \ref{sec:couette} with $\alpha=0$. We will consider multiple cases with different Reynolds numbers, $Re=\{1, 10, 100\}$.
The domain is $\domain = \left[0, 2\pi \right]^2$ with periodic boundary conditions.
We obtain the initial conditions of the conserved state from
\begin{align}
    \rho_0(x,y) = 1 ,\quad 
    p_0(x,y) = p_b - \frac{ M^2 }{4} \left( \cos(2x) + \cos(2y) \right) 
    \nonumber \\
    u_0(x,y) = +M \cos(x)\sin(y) ,\quad
    v_0(x,y) = -M \sin(x)\cos(y) ,
\end{align}
where $M=0.1$ is the Mach number and $p_b = {1}/{\gamma}$ is the base pressure.
At this low Mach number, the flow should evolve similarly to the incompressible limit which we use as our reference solution, 
\begin{equation}
    \rho(x,y,t) = \rho_0(x,y), \qquad
    u(x,y,t) = u_0(x,y) D(t), \qquad 
    v(x,y,t) = v_0(x,y) D(t).
\end{equation}
In the above expression, the decay rate is $D(t) = \exp\left(\frac{-2t}{Re}\right)$.
The total energy is conserved and is equal to its initial value, 
\begin{equation}
    E_T(t) =  \int_{\domain} \rhoe(x,y,0) d \domain = \frac{4\pi^2 p_b}{\gamma-1} + M^2\pi^2
\end{equation}
while the expected total kinetic decay is,
\begin{equation}
    E_K(t) =  \int_{\domain} \frac{1}{2} \rho \left( u_1^2(x,y,t) + u_2^2(x,y,t) \right) d \domain = M^2 \pi^2 D^2(t). 
\end{equation}

We first perform a test with \cEDNN on the whole periodic domain, with periodicity enforced using the feature layer described in \S\ref{sec:boundary} and, hence, without introducing any interface reconstruction.  We then use a \mEDNN approach, where we adopt a $2\times 2$ domain decomposition.  Each sub-domain has four coupled networks for the continuity, the two momentum, and the energy equations.  Each neural network is comprised of four hidden layer with twenty neurons. 
The common solution and normal viscous fluxes were computed with the LDG method and the common inviscid normal fluxes were computed with the Rusanov method (\ref{sec:riemann_solver}).
For both the solution and the flux correction functions, we use $\mathcal{M}_{3,0.2}$. 
Each sub-domain is sampled at $66\times66$ points distributed uniformly, and we march with RK4 and $\Delta t = Re\times10^{-4}$.

\begin{figure}
    \centering
    \begin{subfigure}[b]{0.30\textwidth}
        \centering
        \includegraphics[height=\textwidth]{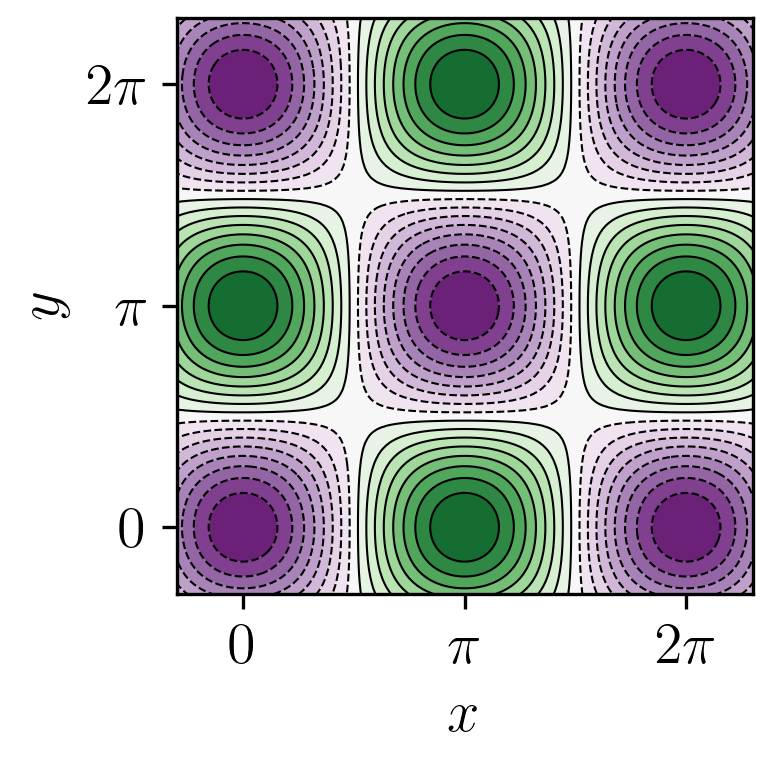}
        \caption{}
    \end{subfigure}
    \begin{subfigure}[b]{.30\textwidth}
        \centering
        \tikz[overlay, remember picture] \node[anchor=south, inner sep=0] (b) 
            {\includegraphics[height=\textwidth]{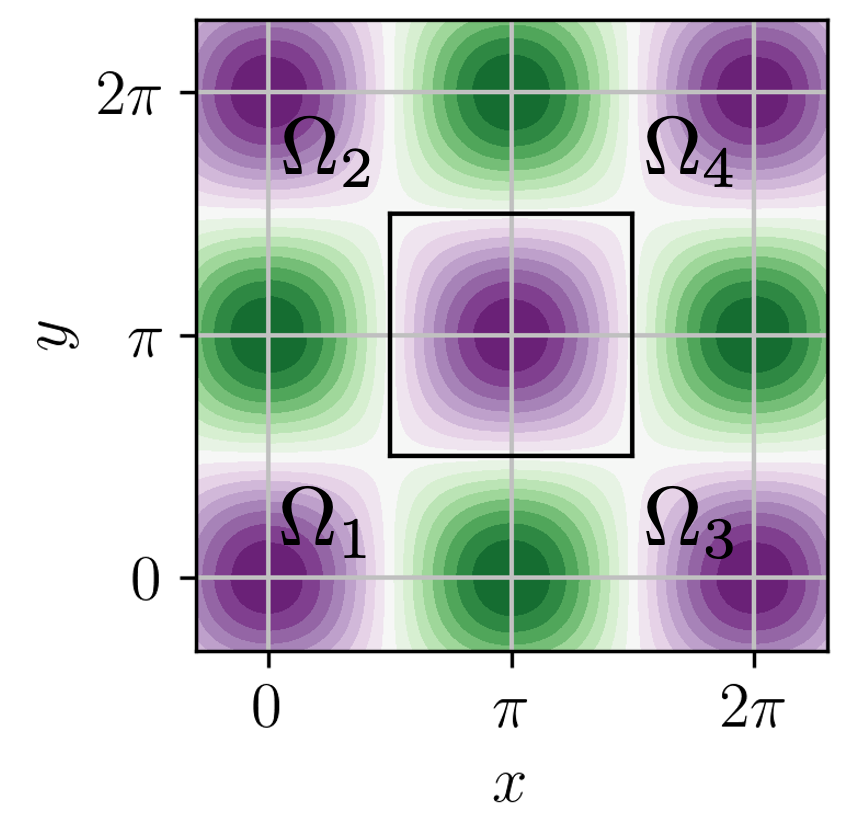}};
        \caption{}
        \label{fig:taylorgreen_2x2}
    \end{subfigure}
    \hspace{0.05\textwidth}
    \begin{subfigure}[b]{.30\textwidth}
        \centering
        \tikz[overlay, remember picture] \node[anchor=south, inner sep=0] (c) 
            {\includegraphics[height=\textwidth]{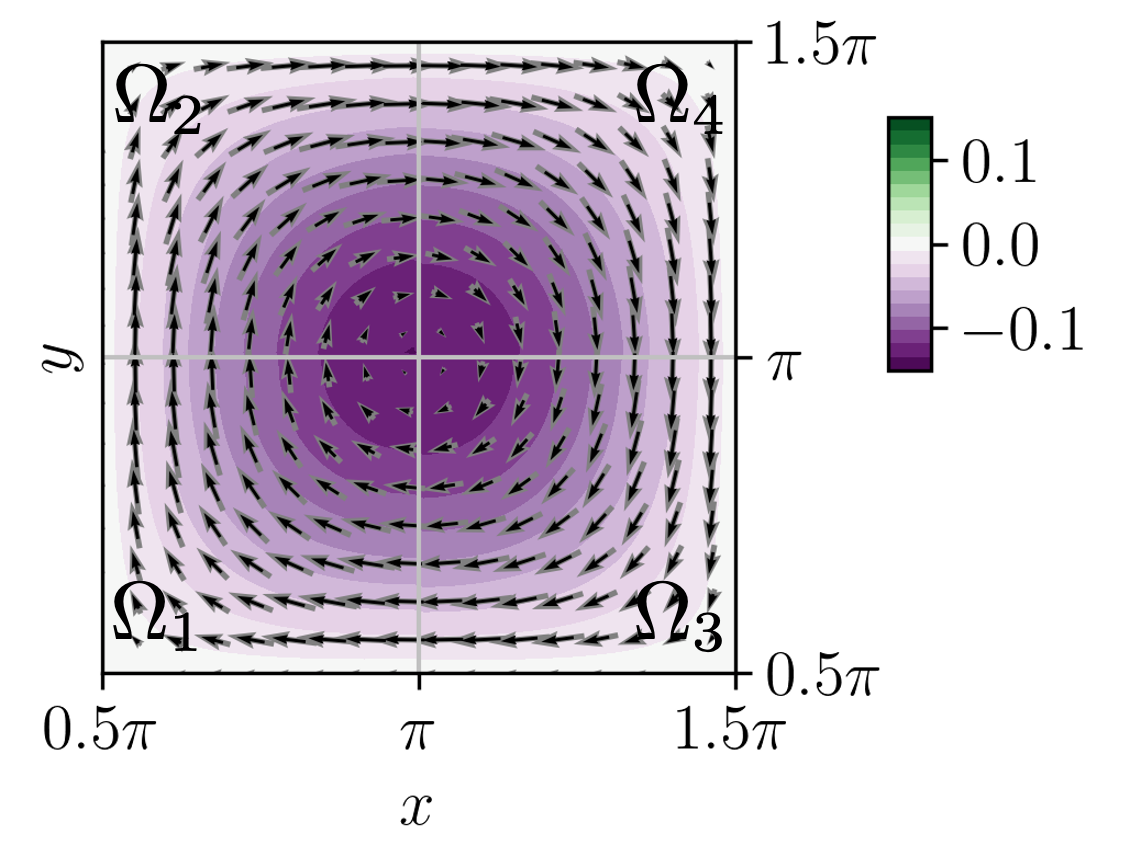}};
        \caption{}
        \label{fig:taylorgreen_zoom}
    \end{subfigure}
    \caption{EDNN prediction of vorticity at $t=20$ for the $Re=100$ case. Color contours show EDNN solution.
    (a) Single \cEDNN with lines showing the reference solution where solid is positive and dashed is negative.
    (b) 2x2 \mEDNN with lines showing the element interfaces.
    (c) Zoom-in on $[0.5\pi, 1.5\pi]^2$ with \mEDNN vectored velocity $[u_1, u_2]$ (black thin) and reference velocity (grey thick).
    }
    \begin{tikzpicture}[overlay, remember picture]
    \draw[black, thin] ([shift={(12.2mm,  11.9mm)}]b.center)--([shift={(6mm, 22.1mm)}]c.west);
    \draw[black, thin] ([shift={(12.2mm, -2.3mm)}]b.center)--([shift={(6mm, -14.1mm)}]c.west);
    \end{tikzpicture}
    \label{fig:taylorgreen_result}
\end{figure}
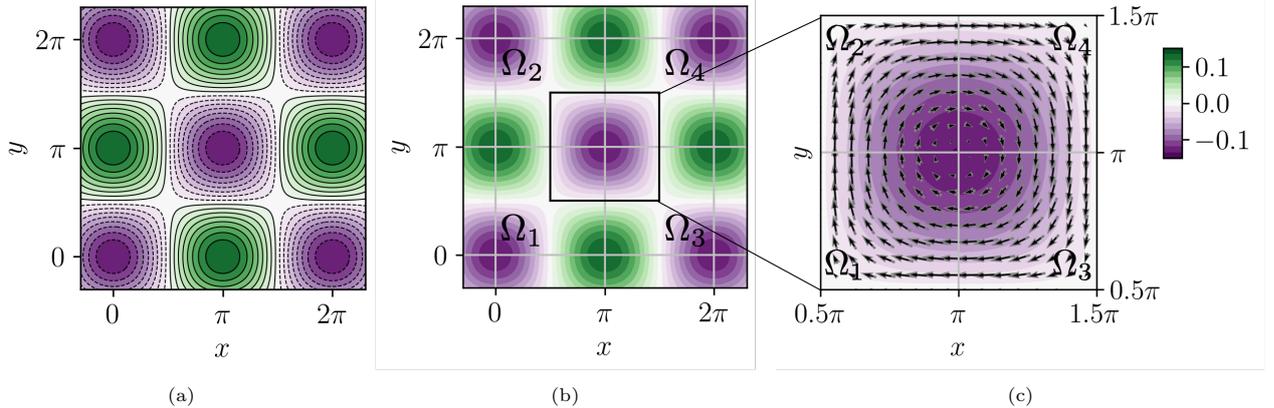

\begin{figure} 
    \centering
    \begin{subfigure}[b]{0.30\textwidth}
        \centering
        \includegraphics[height=\textwidth]{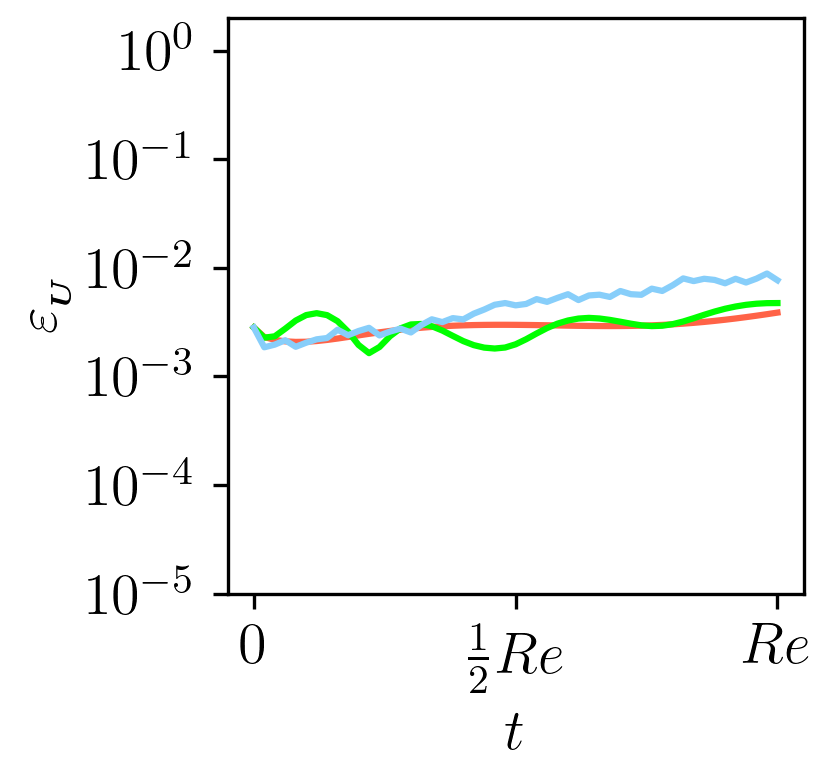}
        \caption{}
        \label{fig:taylorgreen_error}
    \end{subfigure}
    \hspace{0.02\textwidth}
    \begin{subfigure}[b]{0.30\textwidth}
        \centering
        \includegraphics[height=\textwidth]{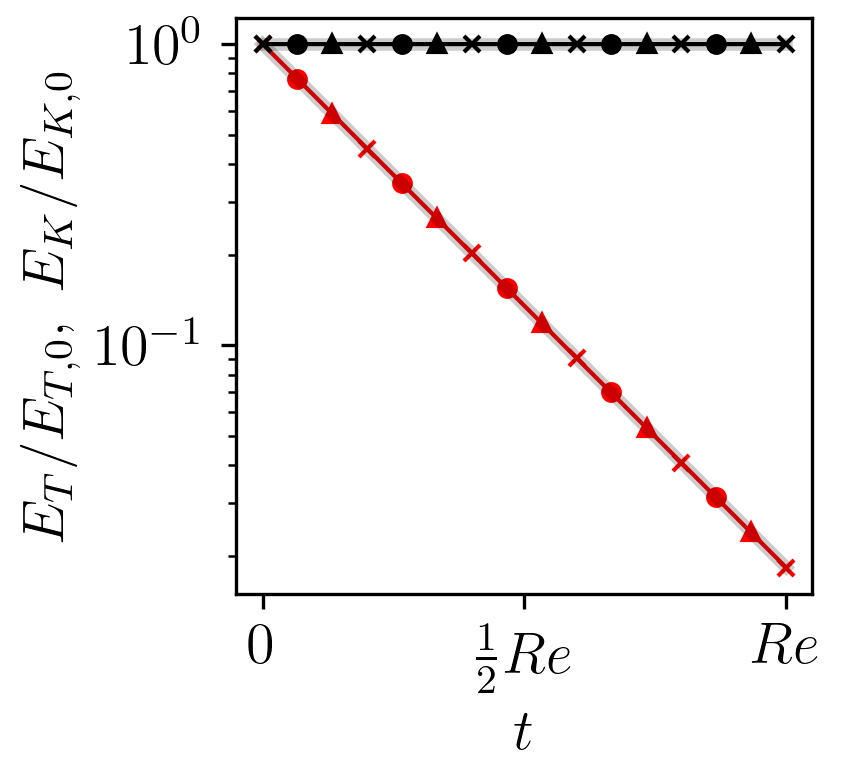}
        \caption{}
        \label{fig:taylorgreen_energy}
    \end{subfigure}
    \caption{Assessment of EDNN error for Taylor-Green vortices. 
    (a) Instantaneous error $\varepsilon$ over time of $\bm{U} = [u_1, u_2]$ velocity for $Re=1$ (red), $Re=10$ (green), and $Re=100$ (blue).
    (b) Change in energy over time of $E_T$ (black) and $E_K$ (red) compared to references (grey thick) for $Re=1$ (circles), $Re=10$ (triangles), and $Re=100$ (crosses).
    }
    \label{}
\end{figure}

In figure \ref{fig:taylorgreen_result}, we report the $Re=100$ case's vorticity, $\omega = \bm{\nabla} \times \bm{u}$, at time $t=20$.  The first panel compares the single \cEDNN solution on the full domain to the analytical expression, using color and line contours respectively.  In panel (b), we report the prediction by the $2\times 2$ \mEDNN with interface corrections both on the internal and also the periodic interfaces between the four sub-domains.  These prediction similarly show very good agreement to the anticipated behaviour of the vorticity.  It is important to note that, while viscous decay is dominant, the nonlinear advection terms in this evolution are all exercised, and must balance exactly or else the vortex pattern is distorted.  The zoomed-in view in figure \ref{fig:taylorgreen_zoom} shows the velocity vectors $\left(u_1, u_2\right)$.  The overlaid black and grey arrows correspond to the \mEDNN prediction and the reference velocity field, and again shows excellent agreement.  

We evaluated the normalized instantaneous error \eqref{eq:error} for the velocity vector $\bm{U} = [u_1, u_2]$.
Figure \ref{fig:taylorgreen_error} shows that the velocities and are accurately predicted.  
The errors at the final times, $t=Re$ for each Reynolds number, are all less than one percent.
In addition, we report the \mEDNN prediction of the total and kinetic energy over time in figure \ref{fig:taylorgreen_energy}, and show that the former remains conserved and the latter decays at the expected rate.

\section{Conclusions}
\label{sec:conclusion}

Use of machine-learning strategies to solve partial differential equations (PDE) is an active area of research.  Evolution equations are of particular interest in the present work, where the solution is expressed by a network in all but one dimension, and evolved in the remaining dimension using the governing equations.  These types of equations are common in physics, where the solution is evolved is  time, although evolution in space for parabolic system can be treated similarly.  The original idea of an evolution deep neural network (EDNN) used a single network to represent the solution of the PDE at an instant in time, and the system of PDEs was then used to determine the time-evolution of the network parameters and, as such, the evolution of the solution for any time horizon of interest.   This approach requires use of a large network that can express the global solution, and is therefore computationally costly.  
In the present work, we use several neural networks whose evolutions are coordinated to accurately predict the solution of the system of PDEs.  First we introduce coupled EDNN (\cEDNN) where we deploy multiple scalar networks each for one of the equations in the system of PDEs, and couple their evolution through the physics coupling in the governing equations.  Second, we introduced distributed EDNN (\dEDNN) where we take advantage of spatial domain decomposition and assign networks to different sub-domains, or elements, in space and coordinate their evolution across interface corrections.
Together, these two strategy form the \mEDNN approach, which is a new framework for numerically solving PDE problems using distributed and coupled neural network.

We demonstrated the capability of \mEDNN for solving canonical PDE problems.   
We used \cEDNN to evolve the compressible Navier-Stokes equations for Couette flow for a range of Mach numbers.  
We showed that \dEDNN successfully applied domain decomposition and flux corrections to solve linear advection and the two dimensional heat equation. 
Finally, \mEDNN coordinated neural networks to solve the Taylor-Green vortices problem.
In all cases, the aggregate solution agreed well with analytical and reference solutions.

Note that, with specific choices of the network architecture, sampling points, and correction functions, EDNN can recover the classical flux reconstruction schemes.
The basic procedure would be to use an input feature layer $x \mapsto [l_{e,1}(x) ,\dots, l_{e,\NUMWEIGHTS}(x) ]$, a single (output) node, a ``linear'' activation function $\sigma(x) = x$, and no bias parameter.  With these choices, the EDNN ansatz becomes 
$\hat{q}_e(x,t) = \sum_{n=1}^\NUMWEIGHTS w_{e,n}(t) l_{e,n}(x)$.
Then consider a fixed set of $\NUMWEIGHTS$ sampling points, $x_{e,n}$, and the Lagrange polynomials, $l_{e,n}$, such that $l_{e,n}(x_{e,n}) = 1$ and zero on the other sampling points.
EDNN is then a degree $\NUMWEIGHTS-1$ polynomial that interpolates the solution at the sampling points.
In this limit, we would not be exploiting the expressivity of the network, and would require a large number of small elements to accurately represent the solution as customary in finite-element methods.  
Instead, by adopting larger networks that have powerful expressivity, we can partition the domain into relatively large elements.    
The combination of network expressively and domain decomposition for distributed networks, which is introduced in the present work, togehter provide the foundation for flexible and computationally efficient solvers.  

In addition, we elected to use generic network architectures and activation functions.  We acknoweldge that it is possible to design networks that are capable of exactly representing, up to machine precision, a known solution of a PDE (e.g.\,using Gaussian activation for an advecting or decaying Gaussian field).  Such networks with problem-specific input feature layers or activation functions can still be adopted with the herein presented approach, and if fact should be used when possible.  However, we assumed no prior expert knowledge of the shape of the solution in our examples.
We demonstrated that \mEDNN remained accurate even while choosing less specialized network architecture, or approximation ansatz.

With \mEDNN,  numerous avenues for future research are worth exploring.
Firstly, \mEDNN can be parallelized to realize computational speedups and to tackle large-scale problems.
Secondly, our approach for imposing boundary conditions is simple, flexible, and can accommodate arbitrary time-dependent geometries.  As such, \mEDNN can be used to simulate flow over immersed bodies, such as canonical flow over a cylinder or a pitching airfoil. 
Thirdly, in the present foundational work, the \mEDNN domain decomposition was limited to conforming hyper-rectagular elements.  The underlying neural networks, however, are not restricted to a subdomain shape and can easily handle irregular geometries. This versatility offers the opportunity for the domain to be partitioned into an unstructured collection of various shapes and sizes.

\section*{Acknowledgements}
\label{sec:acknowledge}
This work was supported by the Defense Advanced Research Projects Agency [grant HR00112220035] and the National Science Foundation [grant 2214925].

\appendix
\section{Adaptive sampling}
\label{sec:adapt_sampling}
At each time-step, EDNN uses a set of sampling points from the spatial domain to solve the linear system \eqref{eq:optimization_linear}.
Electing to use uniformly spaced points is suitable for many problems. 
However, adaptive sampling techniques have been shown to improve the accuracy of EDNN  predictions \cite{Bruna_23, Wen_24, Kast_24}.
We mimic sampling from a probability distribution with a heuristic approach that clusters the points in dynamic regions.
Additionally, we wish to ensure that the global behavior is captured by evolving EDNN with samples that span the entire domain.
To address both, competing goals, our approach approximates a mixture distribution combining the adaptive distribution based on the PDE dynamics with a static, uniform distribution.

\subsection*{Procedure}
Consider a problem in one spatial dimension on $\domain = \left[x_L, x_R \right]$.
Let $x_1^{(t)} <\dots< x_\NUMX^{(t)}$ be the ordered points at time-step $t$.
For the next time-step $t+1$, we seek the positions $x_i^{(t+1)}$ given the data $\adaptHeight_i^{(t)} = \adaptHeight\left(x_i^{(t)}\right)$, where $\adaptHeight$ represents the density function used to determine point locations.

To approximate the area under the curve $\adaptHeight$ at time-step $t$, 
let $\Delta x_i^{(t)} = x_{i+1}^{(t)} - x_i^{(t)}$ be the point spacing, 
$\adaptHeight_i^{(t)} = \frac{1}{2} \left( h_i^{(t)} + h_{i+1}^{(t)} \right)$ be the segment height,
and the total area be,
\begin{equation}
    \adaptArea^{(t)} = \sum_{i=1}^{\NUMX-1} \Delta x_i^{(t)} \adaptHeight_i^{(t)} .
\end{equation}
For the next time step, we wish to find point placements so that all segment areas are equal. 
To avoid implicit solves, we approximate the segment area using the current time-step data, which produces the lengths,
\begin{equation}
    \Delta \hat{x}_i^{(t+1)} 
    = \frac{\adaptArea^{(t)}}{\NUMX-1} \frac{1}{\adaptHeight_i^{(t)}} 
    \approx  \frac{\adaptArea^{(t+1)}}{\NUMX-1} \frac{1}{\adaptHeight_i^{(t+1)}} .
\end{equation}
These segment lengths may not total the domain length, so we normalize,
\begin{equation}
    \Delta x_i^{(t+1)} = \frac{x_R - x_L}{\sum \Delta \hat{x}_j^{(t+1)} }  \Delta \hat{x}_i^{(t+1)}, 
\end{equation}
which gives us the desired point spacing.

In order to approximate the mixture distribution, we choose the data,
\begin{equation}
    \adaptHeight\left(x_i^{(t)}\right) 
    = \adaptParam \Vert\operator(x_i^{(t)}, t, \soln(x_i^{(t)}))\Vert_1 
    + (1-\adaptParam) \frac{1}{\NUMX}\sum_{j=1}^{\NUMX} \Vert\operator(\soln(x_j^{(t)}))\Vert_1,
    \label{eq:adaptive}
\end{equation}
with parameter $\adaptParam \in \left[0,1\right]$.
The first term, $\adaptParam \Vert\operator(x_i^{(t)}, t, \soln(x_i^{(t)}))\Vert_1$, drives the points to the dynamic spatial regions.
The second term drives the points to be uniformly distributed.
The parameter $\adaptParam$ balances these competing goals. 
% Note when $\adaptParam=0$, we recover uniformly spaced points.
For the Couette-flow numerical experiments in \S\ref{sec:couette}, we used $\beta = 0.9.$

This procedure mimics sampling from a distribution $\adaptHeight$ estimated from data $\adaptHeight_i$ without constructing a proper probability distribution estimation, generating random samples, or making additional network inferences.
Also, we can utilize this procedure to determine an initial set of points. 
Using either the initial network prediction or the analytical initial condition, we iteratively adapt the points with the above procedure to settle on the point locations for the first time step.

\section{Riemann solver choices}
\label{sec:riemann_solver}
In this section, we present our choices for computing the common interface solutions and common interface normal flux. 
Overall, we choose upwinding schemes for advective fluxes. 
Since the viscous flux uses the corrected solution, we pair the common interface solution and the common viscous flux solvers so that the viscous flux remains centered.
To simplify, we adopt the following notation for the local (subscript minus) and external (subscript plus) fluxes,
\begin{equation}
    \begin{split}
        \flux_{inv,-} = \fluxfunc_{inv} (\soln_{-}) ,\qquad &
        \flux_{vis,-} = \fluxfunc_{vis} (\soln_{-}, \aux_{-}) , \\
        \flux_{inv,+} = \fluxfunc_{inv} (\soln_{+}) ,\qquad &
        \flux_{vis,+} = \fluxfunc_{vis} (\soln_{+}, \aux_{+}) .
    \end{split}
\end{equation}

\subsection*{Linear advection and diffusion}
For the common interface solution, we use the local discontinuous Galerkin (LDG) \cite{Cockburn_LDG} (upwinding) method:
\begin{equation}
    \solnSolver (\soln_{-}, \soln_{+}, \bm{n}) = 
    \frac{1}{2} (\soln_{-} + \soln_{+}) 
    +  \, \text{sgn}\left(\bm{u} \cdot \bm{n}\right) \frac{1}{2} (\soln_{-} - \soln_{+}),
\end{equation}
where $\bm{u}$ is the problem wave direction and $\text{sgn}$ is the sign function,
\begin{equation}
    \text{sgn}(x) = 
    \begin{cases}
     1 &\text{if } x > 0 \\
     0 &\text{if } x = 0 \\
    -1 &\text{if } x < 0 \\
    \end{cases} .
\end{equation}
Note that, if $\bm{u} \cdot \bm{n} > 0 $ then $\soln^\star = \soln_{-}^\dis$,
and if  $\bm{u} \cdot \bm{n} < 0 $ then $\soln^\star = \soln_{+}^\dis$ .

For the common inviscid normal flux, we use the Lax-Friedrichs method \cite{Friedrichs_LF}:
\begin{equation}
    \fluxSolver_{inv}(\soln_{-}, \soln_{+}, \bm{n})  = 
    \frac{1}{2}\left( \flux_{inv,-} + \flux_{inv,+} \right) \cdot  \bm{n}
    + \frac{r_1}{2} |\bm{u} \cdot \bm{n}| \left( \soln_{-} - \soln_{+} \right),
\end{equation}
where $r_1 \in [0,1]$ is a parameter, $r_1=0$ recovers an central average flux, and $r_1=1$ recovers a fully upwind flux.
In our experiments, we choose $r_1 = 1$.

For the common viscous normal flux, we use the local discontinuous Galerkin (LDG) (downwinding) method:
\begin{equation}
    \fluxSolver_{vis}(\soln_{-}, \soln_{+}, \aux_{-}, \aux_{+}, \bm{n}) = 
    \frac{1}{2}\left( \flux_{vis,-} + \flux_{vis,+} \right) \cdot \bm{n}
    - \text{sgn}\left(\bm{u} \cdot \bm{n}\right) \frac{1}{2}\left(\flux_{vis,-} - \flux_{vis,+} \right) \cdot \bm{n} 
    + r_2 \left( \soln_{-} - \soln_{+} \right),
\end{equation}
where $r_2$ is a parameter penalizing the jump in the solution.
Note that the use of a downwinding method for the viscous fluxes is to compensate for the upwinding performed on the common interface solution. 
This pairing results in a central viscous flux scheme.
In our experiments, we choose $r_2 = 0.1$.

\subsection*{Compressible Navier-Stokes}
For the common inviscid normal flux, we use the Rusanov method \cite{Rusanov} with estimated wavespeed from \citet{Witherden_pyfr}:
\begin{equation}
    \fluxSolver_{inv}(\soln_{-}, \soln_{+}, \bm{n})  = 
    \frac{1}{2} \left( \flux_{inv,-} + \flux_{inv,+} \right) \cdot \bm{n}
    + \frac{s}{2} \left( \soln_{-} - \soln_{+} \right),
\end{equation}
where $s$ is an estimate of the maximum wave speed,
\begin{equation}
    s 
    = \sqrt{\frac{\gamma \left(p_{-} + p_{+} \right)}{\rho_{-} + \rho_{+}} }
    + \frac{1}{2} \left| \bm{n} \cdot \left( \bm{u}_{-} + \bm{u}_{+} \right) \right|,
\end{equation}
and $\rho_{\mp}$ are the local and external densities, 
$p_{\mp}$ are the local and external pressures, 
$\bm{u}_{\mp}$ are the local and external velocity vectors,
and $\gamma$ is the ratio of specific heats. 

For the common solution and common viscous normal flux, we use the LDG methods.
However, instead of a prescribed velocity, we use the predicted, pointwise velocity averaged from the local and external boundary state,
$\bm{u} = 0.5 \left(\bm{u}_{-}+\bm{u}_{+}\right)$.

\section{Evaluation of Navier-Stokes viscous flux}
\label{sec:cns_viscous}
In this section we provide an example of the computation of the divergence of flux for \mEDNN. 
The viscous flux \eqref{eq:flux} is a function of the discontinuous conserved state and the (corrected) auxiliary variable.
Therefore divergence of the viscous flux \eqref{eq:divergence_vis_flux} requires the gradient of the discontinuous state $\nabla \bm{q}^D$ and the auxiliary variable $\bm{a}^D$.
Though both quantities are approximations of the solution gradient, they are not equal.
To highlight this difference, we mark  variables that depend on the discontinuous state $\left( \soln^\dis, \nabla\soln^\dis \right)$ alone with overline $\overline{\Box}$ and name them `discontinuous'.
In contrast, we indicate variables that additionally depend on the auxiliary $\left( \soln^\dis, \nabla\soln^\dis, \aux^\dis, \nabla\aux^\dis \right)$ with tilde
$\widetilde{\Box}$ and name them `corrected'. 
Certain physical variables will have both a discontinuous and corrected versions.

We rewrite the Navier-Stokes equations \eqref{eq:cns_flux} with this notation:
\begin{equation}
    \ddfrac{t}{\soln} + {\nabla\cdot} \overline{\fluxfunc}_{inv} + {\nabla\cdot} \widetilde{\fluxfunc}_{vis} = 0
\end{equation}
\begin{equation}
    \overline{\fluxfunc}_{inv,i}
    =
    \begin{bmatrix}
    \overline{\rho u}_i  \\
    \overline{\rhou}_s \overline{u}_i + \overline{p}\delta_{si} \\
    \overline{u}_i \left( \overline{\rhoe} + \overline{p} \right) \\
    \end{bmatrix}
    , \qquad
    \widetilde{\fluxfunc}_{vis,i}
    =
    \begin{bmatrix}
    0 \\
    - \widetilde{\tau}_{si} \\
    - \overline{u}_k \widetilde{\tau}_{ki} + \widetilde{\theta}_i \\
    \end{bmatrix} 
    ,
\end{equation}
\begin{equation}
    \overline{u}_i = \frac{\overline{\rhou}_i }{\overline{\rho}} 
    ,\qquad
    \overline{p} = \left( \gamma-1 \right) \left( \overline{\rhoe} - \frac{1}{2} \overline{\rhou}_k \, \overline{u}_k\right) 
    ,\qquad
    \overline{T} = \frac{\overline{p}}{R_g \overline{\rho}}
    ,\qquad
    \overline{\mu} = \left( \left( \gamma -1 \right) \overline{T} \right)^\alpha ,
\end{equation}
\begin{equation}
    \widetilde{\tau}_{ki}
    = \frac{\overline{\mu}}{Re} \left( \ddfrac{x_k}{\widetilde{u}_i} + \ddfrac{x_i}{\widetilde{u}_k} \right) 
    + \frac{1}{Re}\left(\mu_b - \frac{S-1}{S }\overline{\mu} \right) \ddfrac{x_l}{\widetilde{u}_l} \delta_{ki}
    , \qquad
    \widetilde{\theta}_i = - \frac{\overline{\mu}}{Re Pr} \ddfrac{x_i}{\widetilde{T}} .
\end{equation}
Consider, for example, the divergence of viscous energy flux,
\begin{equation}
    \nabla \cdot \widetilde{\fluxfunc}_{vis,\rhoe} 
    = 
    - \widetilde{\tau_{ki}} \ddfrac{x_i}{\overline{u_k}}
    - \overline{u_k} \ddfrac{x_i}{\widetilde{\tau_{ki}}}
    + \ddfrac{x_i}{\widetilde{\theta_i}} .
\end{equation}
Notice the discontinuous velocity gradient appears, and from the viscous stress tensor, the corrected velocity gradient is needed.
The discontinuous velocity gradient,
\begin{equation}
    \ddfrac{x_j}{\overline{u_i}}
    = - \overline{\rho}^{-2} \overline{\rhou_i} \ddfrac{x_j}{\overline{\rho}}
      + \overline{\rho}^{-1} \ddfrac{x_j}{\overline{\rhou_i}},
\end{equation}
can be straightforwardly computed from the discontinuous state using the EDNN solution and automatic differentiation.
However, the corrected velocity gradient,
\begin{equation}
    \ddfrac{x_j}{\widetilde{u_i}}
    = - \overline{\rho}^{-2} \overline{\rhou_i} \ddfrac{x_j}{\widetilde{\rho}}
      + \overline{\rho}^{-1}                    \ddfrac{x_j}{\widetilde{\rhou_i}},
\end{equation}
requires the auxiliary of density and momentum.
The computation of all other divergences and their necessary components follows similarly.

\bibliographystyle{unsrtnat} 
\bibliography{references}

\begin{thebibliography}{29}
\providecommand{\natexlab}[1]{#1}
\providecommand{\url}[1]{\texttt{#1}}
\expandafter\ifx\csname urlstyle\endcsname\relax
  \providecommand{\doi}[1]{doi: #1}\else
  \providecommand{\doi}{doi: \begingroup \urlstyle{rm}\Url}\fi

\bibitem[Raissi et~al.(2019)Raissi, Perdikaris, and Karniadakis]{Raissi_19}
Maziar Raissi, Paris Perdikaris, and George~E Karniadakis.
\newblock Physics-informed neural networks: A deep learning framework for solving forward and inverse problems involving nonlinear partial differential equations.
\newblock \emph{Journal of Computational physics}, 378:\penalty0 686--707, 2019.
\newblock \doi{https://doi.org/10.1016/j.jcp.2018.10.045}.
\newblock URL \url{https://www.sciencedirect.com/science/article/pii/S0021999118307125}.

\bibitem[Jagtap and Karniadakis(2020)]{Jagtap_20}
Ameya~D Jagtap and George~Em Karniadakis.
\newblock Extended physics-informed neural networks (xpinns): A generalized space-time domain decomposition based deep learning framework for nonlinear partial differential equations.
\newblock \emph{Communications in Computational Physics}, 28\penalty0 (5), 2020.
\newblock \doi{https://doi.org/10.4208/cicp.oa-2020-0164}.
\newblock URL \url{https://www.osti.gov/biblio/2282003}.

\bibitem[Clark Di~Leoni et~al.(2023)Clark Di~Leoni, Agarwal, Zaki, Meneveau, and Katz]{Patricio_23}
Patricio Clark Di~Leoni, Karuna Agarwal, Tamer~A. Zaki, Charles Meneveau, and Joseph Katz.
\newblock Reconstructing turbulent velocity and pressure fields from under-resolved noisy particle tracks using physics-informed neural networks.
\newblock \emph{Experiments in Fluids}, 64\penalty0 (5):\penalty0 95, 2023.
\newblock \doi{10.1007/s00348-023-03629-4}.

\bibitem[Lu et~al.(2021)Lu, Jin, Pang, Zhang, and Karniadakis]{Lu_21}
Lu~Lu, Pengzhan Jin, Guofei Pang, Zhongqiang Zhang, and George~Em Karniadakis.
\newblock Learning nonlinear operators via deeponet based on the universal approximation theorem of operators.
\newblock \emph{Nature machine intelligence}, 3\penalty0 (3):\penalty0 218--229, 2021.

\bibitem[Li et~al.(2020)Li, Kovachki, Azizzadenesheli, Bhattacharya, Stuart, Anandkumar, et~al.]{Li_20}
Zongyi Li, Nikola~Borislavov Kovachki, Kamyar Azizzadenesheli, Kaushik Bhattacharya, Andrew Stuart, Anima Anandkumar, et~al.
\newblock Fourier neural operator for parametric partial differential equations.
\newblock In \emph{International Conference on Learning Representations}, 2020.

\bibitem[Cai et~al.(2021)Cai, Wang, Lu, Zaki, and Karniadakis]{Cai_21}
S.~Cai, Z.~Wang, L.~Lu, T.~A. Zaki, and G.~E. Karniadakis.
\newblock {DeepM\&Mnet}: Inferring the electroconvection multiphysics fields based on operator approximation by neural networks.
\newblock \emph{Journal of Computational Physics}, 436:\penalty0 110296, 2021.
\newblock ISSN 0021-9991.
\newblock \doi{https://doi.org/10.1016/j.jcp.2021.110296}.
\newblock URL \url{https://www.sciencedirect.com/science/article/pii/S0021999121001911}.

\bibitem[Du and Zaki(2021)]{Du_ednn}
Yifan Du and Tamer~A. Zaki.
\newblock Evolutional deep neural network.
\newblock \emph{Phys. Rev. E}, 104:\penalty0 045303, Oct 2021.
\newblock \doi{10.1103/PhysRevE.104.045303}.
\newblock URL \url{https://link.aps.org/doi/10.1103/PhysRevE.104.045303}.

\bibitem[Anderson and Farazmand(2024)]{Anderson_24}
William Anderson and Mohammad Farazmand.
\newblock Fast and scalable computation of shape-morphing nonlinear solutions with application to evolutional neural networks.
\newblock \emph{Journal of Computational Physics}, 498:\penalty0 112649, 2024.
\newblock ISSN 0021-9991.
\newblock \doi{https://doi.org/10.1016/j.jcp.2023.112649}.
\newblock URL \url{https://www.sciencedirect.com/science/article/pii/S0021999123007441}.

\bibitem[Bruna et~al.(2024)Bruna, Peherstorfer, and Vanden-Eijnden]{Bruna_23}
Joan Bruna, Benjamin Peherstorfer, and Eric Vanden-Eijnden.
\newblock Neural galerkin schemes with active learning for high-dimensional evolution equations.
\newblock \emph{Journal of Computational Physics}, 496:\penalty0 112588, 2024.
\newblock ISSN 0021-9991.
\newblock \doi{https://doi.org/10.1016/j.jcp.2023.112588}.
\newblock URL \url{https://www.sciencedirect.com/science/article/pii/S0021999123006836}.

\bibitem[Wen et~al.(2024)Wen, Vanden-Eijnden, and Peherstorfer]{Wen_24}
Yuxiao Wen, Eric Vanden-Eijnden, and Benjamin Peherstorfer.
\newblock Coupling parameter and particle dynamics for adaptive sampling in neural galerkin schemes.
\newblock \emph{Physica D: Nonlinear Phenomena}, 462:\penalty0 134129, 2024.
\newblock ISSN 0167-2789.
\newblock \doi{https://doi.org/10.1016/j.physd.2024.134129}.
\newblock URL \url{https://www.sciencedirect.com/science/article/pii/S0167278924000800}.

\bibitem[Kast and Hesthaven(2024)]{Kast_24}
Mariella Kast and Jan~S. Hesthaven.
\newblock Positional embeddings for solving pdes with evolutional deep neural networks.
\newblock \emph{Journal of Computational Physics}, 508:\penalty0 112986, 2024.
\newblock ISSN 0021-9991.
\newblock \doi{https://doi.org/10.1016/j.jcp.2024.112986}.
\newblock URL \url{https://www.sciencedirect.com/science/article/pii/S0021999124002353}.

\bibitem[Chen et~al.(2023)Chen, Wu, Grinspun, Zheng, and Chen]{Chen_23}
Honglin Chen, Rundi Wu, Eitan Grinspun, Changxi Zheng, and Peter~Yichen Chen.
\newblock Implicit neural spatial representations for time-dependent {PDE}s.
\newblock In Andreas Krause, Emma Brunskill, Kyunghyun Cho, Barbara Engelhardt, Sivan Sabato, and Jonathan Scarlett, editors, \emph{Proceedings of the 40th International Conference on Machine Learning}, volume 202 of \emph{Proceedings of Machine Learning Research}, pages 5162--5177. PMLR, 23--29 Jul 2023.
\newblock URL \url{https://proceedings.mlr.press/v202/chen23af.html}.

\bibitem[Finzi et~al.(2023)Finzi, Potapczynski, Choptuik, and Wilson]{Finzi_23}
Marc~Anton Finzi, Andres Potapczynski, Matthew Choptuik, and Andrew~Gordon Wilson.
\newblock A stable and scalable method for solving initial value {PDE}s with neural networks.
\newblock In \emph{The Eleventh International Conference on Learning Representations}, 2023.
\newblock URL \url{https://openreview.net/forum?id=vsMyHUq_C1c}.

\bibitem[Kao et~al.(2024)Kao, Zhao, and Zhang]{Kao_24}
Tunan Kao, Jin Zhao, and Lei Zhang.
\newblock petnns: Partial evolutionary tensor neural networks for solving time-dependent partial differential equations.
\newblock \emph{arXiv preprint arXiv:2403.06084}, 2024.

\bibitem[Ben-Nun and Hoefler(2019)]{Ben_19}
Tal Ben-Nun and Torsten Hoefler.
\newblock Demystifying parallel and distributed deep learning: An in-depth concurrency analysis.
\newblock \emph{ACM Comput. Surv.}, 52\penalty0 (4), aug 2019.
\newblock ISSN 0360-0300.
\newblock \doi{10.1145/3320060}.
\newblock URL \url{https://doi.org/10.1145/3320060}.

\bibitem[Merzari et~al.(2023)Merzari, Hamilton, Evans, Min, Fischer, Kerkemeier, Fang, Romano, Lan, Phillips, et~al.]{Merzari_23}
Elia Merzari, Steven Hamilton, Thomas Evans, Misun Min, Paul Fischer, Stefan Kerkemeier, Jun Fang, Paul Romano, Yu-Hsiang Lan, Malachi Phillips, et~al.
\newblock Exascale multiphysics nuclear reactor simulations for advanced designs.
\newblock In \emph{Proceedings of the International Conference for High Performance Computing, Networking, Storage and Analysis}, pages 1--11, 2023.

\bibitem[Becker et~al.(1981)Becker, Carey, and Oden]{becker_fem}
E.B. Becker, G.F. Carey, and J.T. Oden.
\newblock \emph{Finite Elements: An introduction}.
\newblock Number v. 1 in ACM monograph series. Prentice-Hall, 1981.
\newblock ISBN 9780133170573.
\newblock URL \url{https://books.google.com/books?id=Qh3BugEACAAJ}.

\bibitem[Cockburn et~al.(2012)Cockburn, Karniadakis, and Shu]{Cockburn_DG}
Bernardo Cockburn, George~E Karniadakis, and Chi-Wang Shu.
\newblock \emph{Discontinuous Galerkin methods: theory, computation and applications}, volume~11.
\newblock Springer Science \& Business Media, 2012.

\bibitem[Huynh(2007)]{Huynh_fr}
H.~T. Huynh.
\newblock A flux reconstruction approach to high-order schemes including discontinuous galerkin methods.
\newblock In \emph{18th AIAA Computational Fluid Dynamics Conference}. American Institute of Aeronautics and Astronautics, 2007.
\newblock \doi{10.2514/6.2007-4079}.
\newblock URL \url{https://arc.aiaa.org/doi/abs/10.2514/6.2007-4079}.

\bibitem[Vincent et~al.(2011)Vincent, Castonguay, and Jameson]{Vincent_fr}
Peter~E Vincent, Patrice Castonguay, and Antony Jameson.
\newblock A new class of high-order energy stable flux reconstruction schemes.
\newblock \emph{Journal of Scientific Computing}, 47:\penalty0 50--72, 2011.

\bibitem[Castonguay et~al.(2013)Castonguay, Williams, Vincent, and Jameson]{Castonguay_13}
Patrice Castonguay, David~M Williams, Peter~E Vincent, and Antony Jameson.
\newblock Energy stable flux reconstruction schemes for advection--diffusion problems.
\newblock \emph{Computer Methods in Applied Mechanics and Engineering}, 267:\penalty0 400--417, 2013.

\bibitem[Witherden et~al.(2014)Witherden, Farrington, and Vincent]{Witherden_pyfr}
Freddie~D Witherden, Antony~M Farrington, and Peter~E Vincent.
\newblock Pyfr: An open source framework for solving advection--diffusion type problems on streaming architectures using the flux reconstruction approach.
\newblock \emph{Computer Physics Communications}, 185\penalty0 (11):\penalty0 3028--3040, 2014.

\bibitem[Bj{\"o}rck(1996)]{bjorck}
{\AA}ke Bj{\"o}rck.
\newblock \emph{Numerical methods for least squares problems}.
\newblock SIAM, 1996.

\bibitem[Yazdani et~al.(2020)Yazdani, Lu, Raissi, and Karniadakis]{Yazdani}
Alireza Yazdani, Lu~Lu, Maziar Raissi, and George~Em Karniadakis.
\newblock Systems biology informed deep learning for inferring parameters and hidden dynamics.
\newblock \emph{PLOS Computational Biology}, 16\penalty0 (11):\penalty0 1--19, 11 2020.
\newblock \doi{10.1371/journal.pcbi.1007575}.
\newblock URL \url{https://doi.org/10.1371/journal.pcbi.1007575}.

\bibitem[Bodony(2006)]{Bodony}
Daniel~J Bodony.
\newblock Analysis of sponge zones for computational fluid mechanics.
\newblock \emph{Journal of Computational Physics}, 212\penalty0 (2):\penalty0 681--702, 2006.

\bibitem[Anderson and Farazmand(2022)]{Anderson_22}
William Anderson and Mohammad Farazmand.
\newblock Evolution of nonlinear reduced-order solutions for pdes with conserved quantities.
\newblock \emph{SIAM Journal on Scientific Computing}, 44\penalty0 (1):\penalty0 A176--A197, 2022.
\newblock \doi{10.1137/21M1415972}.
\newblock URL \url{https://doi.org/10.1137/21M1415972}.

\bibitem[Cockburn and Shu(1998)]{Cockburn_LDG}
Bernardo Cockburn and Chi-Wang Shu.
\newblock The local discontinuous galerkin method for time-dependent convection-diffusion systems.
\newblock \emph{SIAM journal on numerical analysis}, 35\penalty0 (6):\penalty0 2440--2463, 1998.

\bibitem[Friedrichs(1954)]{Friedrichs_LF}
Kurt~O Friedrichs.
\newblock Symmetric hyperbolic linear differential equations.
\newblock \emph{Communications on pure and applied Mathematics}, 7\penalty0 (2):\penalty0 345--392, 1954.

\bibitem[Rusanov(1962)]{Rusanov}
Vladimir~Vasil'evich Rusanov.
\newblock \emph{Calculation of interaction of non-steady shock waves with obstacles}.
\newblock NRC, Division of Mechanical Engineering, 1962.

\end{thebibliography}

\end{document}